\documentclass[11pt]{jsg}
\usepackage{amssymb,epic,eepic,verbatim,graphicx}
\renewcommand\theenumi{\alph{enumi}}
\renewcommand\theenumii{\roman{enumii}}
\newtheorem{theorem}{Theorem}[subsection]

\newtheorem{corollary}[theorem]{Corollary}

\newtheorem{proposition}[theorem]{Proposition}
\newtheorem{lemma}[theorem]{Lemma}
\newtheorem{lem}[theorem]{}
\theoremstyle{definition}
\newtheorem{definition}[theorem]{Definition}
\theoremstyle{remark}

\newtheorem{example}[theorem]{Example}
\newcommand{\blem}{\begin{lem} \rm}
\newcommand{\elem}{\end{lem}}
\newcommand\A{\mathcal{A}}

\newcommand\M{\mathcal{M}}
\newcommand\D{\mathcal{D}}

\renewcommand\M{\mathcal{M}}

\renewcommand\S{\mathcal{S}}

\renewcommand{\L}{\mathcal{L}}
\newcommand{\I}{\mathcal{I}}
\renewcommand{\O}{\mathcal{O}}

\newcommand{\F}{\mathcal{F}}

\newcommand{\R}{\mathbb{R}}
\renewcommand{\H}{\mathcal{H}}
\newcommand{\RR}{\mathcal{R}}
\newcommand{\C}{\mathbb{C}}

\newcommand{\Z}{\mathbb{Z}}

\newcommand{\ddt}{\frac{d}{dt}}

\renewcommand{\P}{\mathbb{P}}

\newcommand\lie[1]{\mathfrak{#1}}
\renewcommand{\k}{\lie{k}}

\newcommand{\g}{\lie{g}}

\newcommand{\p}{\lie{p}}

\newcommand{\z}{\lie{z}}
\renewcommand{\t}{\lie{t}}

\newcommand{\so}{\lie{so}}
\renewcommand{\sl}{\lie{sl}}
\newcommand{\on}{\operatorname}

\newcommand{\Fr}{\on{Fr}}
\newcommand{\Sp}{\on{Sp}}
\newcommand{\sing}{\on{sing}}
\newcommand{\reg}{\on{reg}}
\newcommand{\Pf}{\on{Pf}}

\newcommand{\height}{\on{ht}}
\newcommand{\Aut}{ \on{Aut} }

\newcommand{\Hom}{ \on{Hom}}
\newcommand{\Ind}{ \on{Ind}}

\newcommand{\Res}{ \on{Res}}

\newcommand{\Vol}{  \on{Vol}}
\newcommand{\RVol}{  \on{RVol}}

\newcommand{\codim}{\on{codim}}

\newcommand{\ssm}{\kern-.5ex \smallsetminus \kern-.5ex}

\newcommand\dirac{/\kern-1.2ex\partial} 
\newcommand\qu{/\kern-.7ex/} 
\newcommand\lqu{\backslash \kern-.7ex \backslash} 

\newcommand\dr{r_+ \kern-.7ex - \kern-.7ex r_-}
 
\newcommand{\labell}\label

\renewcommand{\d}{{\mbox{d}}}
\newcommand{\ol}{\overline}

\newcommand\Phinv{\Phi^{-1}}

\newcommand\eps{\epsilon}

\newcommand\om{\omega}

\newcommand{\f}{\frac}

\newcommand{\hh}{{\f{1}{2}}}

\newcommand{\qq}{{\f{1}{4}}}

\newcommand{\ti}{\tilde}

\renewcommand\vol{\on{vol}}

\newcommand\curv{\on{curv}}

\newcommand\Eul{\on{Eul}}

\newcommand\Vect{\on{Vect}}

\newcommand\grad{\on{grad}}

\newcommand\hull{\on{hull}}
\newcommand\crit{\on{crit}}

\newcommand\bdefn{\begin{definition}}
\newcommand\edefn{\end{definition}}
\newcommand\bea{\begin{eqnarray*}}
\newcommand\eea{\end{eqnarray*}}
\newcommand\bcv{\left[ \begin{array}{r} }
\newcommand\ecv{\end{array} \right] }

\newcommand\bma{\left[ \begin{array} }
\newcommand\ema{\end{array} \right]}
\newcommand\bsj{\left\{ \begin{array}{rrr} }
\newcommand\esj{\end{array} \right\}}

\newcommand\ben{\begin{enumerate}}
\newcommand\een{\end{enumerate}}
\newcommand\bex{\begin{example}}
\newcommand\eex{$\Diamond$ \end{example}}

\begin{document}

\title[Localization for the norm-square of the moment map]{Localization for
the norm-square of the moment map and the two-dimensional Yang-Mills
integral}

\author{Chris T. Woodward} \thanks{Partially supported by NSF
 grant DMS0093647}

\address{Mathematics-Hill Center,
Rutgers University, 110 Frelinghuysen Road, Piscataway, NJ 08854-8019,
U.S.A.}  \email{ctw@math.rutgers.edu}

\date{May 12, 2005}

\maketitle


\section{Introduction}
\label{introsec}

The first seven sections of the paper contain a version of
localization for the norm-square of the moment map in equivariant de
Rham theory.  The main Theorem \ref{Hform} expresses the push-forward
of an equivariant cohomology class on a Hamiltonian $K$-manifold with
proper moment map as a sum of contributions from fixed point
components of one-parameter subgroups corresponding to the critical
values of the norm-square of the moment map.  If the critical set of
the norm-square is non-degenerate in a sense explained below, there is
an improved result Theorem \ref{Hsmooth} which expresses the
push-forward as an integral over the quotient of the critical set.
P.-E. Paradan obtained the same results but did not entirely publish;
most of the ingredients as well as similar results appear in his
papers \cite{pa:loc}, \cite{pa:mo}, and \cite{pa:lo}.  The proof given
here is different from Paradan's.  The existence of a localization
formula, but not the precise form of the contributions, was first
suggested by Witten in his study \cite{wi:tw} of two-dimensional
Yang-Mills theory.  Later Jeffrey and Kirwan \cite{je:lo1} gave a
formula which had a similar purpose but was expressed in terms of
rather different fixed point data.  A $K$-theory version was given by
Vergne \cite{ve:mu} and Paradan \cite{pa:lo}.  A similar result for
sheaf cohomology that I learned from C. Teleman is explained in the
eighth section.

The ninth section contains a definition and computation of the
Yang-Mills path integral in two dimensions.  The idea is to reverse
the logic in Witten's paper \cite{wi:tw}, and take the ``stationary
phase approximation'' (that is, the localization formula) as the
definition of the path integral.  In order to compute it we apply a
symmetry argument from Teleman-Woodward \cite{te:in} which reduces the
computation to an integration over Jacobians.  The result is what the
physicists call the Migdal formula for the path integral; its large
coupling (topological) limit is the Witten volume formula.  More
general formulas for intersection pairings and indices on the moduli
space are given in \cite{te:in}, Meinrenken \cite{me:onwi}, and
Jeffrey-Kirwan \cite{je:in}.  A ``combinatorial'' definition and
computation of the 2d Yang-Mills measure, including observables, is
given by Levy \cite{le:ym}.  Putting these results together shows that
the ``stationary phase'' and ``combinatorial'' definitions of the
two-dimensional Yang-Mills integral are equal.  This might be seen as
the two-dimensional analog of the much harder conjecture regarding the
three-dimensional Chern-Simons path integral, that the
``combinatorial'' definition via Reshetikhin-Turaev \cite{rt:in}
agrees with the ``stationary phase'' definition of Axelrod-Singer
\cite{as:cs2} (with leading order term given as in \cite{freed:comp}).

An appendix contains a proof of an unpublished result of Duistermaat
that the gradient flow of minus the norm-square of the moment map
converges.  

{\em Acknowledgments.} Thanks to the Mathematics Department of the
University of Otago, Dunedin, New Zealand for its hospitality while
writing this paper, to P.-E. Paradan for explaining his work to me and
pointing out a number of mistakes in an earlier version, and to
M. Harada, E. Meinrenken, and C. Teleman for helpful discussions.

\section{Localization for one-parameter subgroups}
\label{cohsec}

Let $K$ be a compact connected Lie group with Lie algebra $\k$ and $M$
a $K$-manifold.  The equivariant de Rham cohomology $H_K(M)$ of $M$
(complex coefficients) can be computed in the {\em Cartan model},
$$ \Omega_K(M) := (S(\k^*) \otimes \Omega(M))^K $$
where $\Omega(M)$ denotes the space of smooth forms on $M$ and
$S(\k^*)$ the symmetric algebra on $\k^*$, see \cite{gu:eqdr}.  The
equivariant differential $d_K$ can be written
$$ d_K: \ \Omega_K(M) \to \Omega_K(M) , \ \ \ (d_K \eta)(\zeta) =
(d + 2 \pi i \iota(\zeta_M))(\eta(\zeta)), \ \ \ \zeta \in \k $$
where $\zeta_M \in \Vect(M)$ denotes the generating vector field
$\zeta_M(m) = [\exp( - t \zeta) m ] $ and $\iota(\zeta_M)$ contraction
with $\zeta_M$.

Suppose $K$ acts locally freely on $M$; then the equivariant
cohomology $H_K(M)$ is isomorphic to $H(K \backslash M)$ via pullback
$p^*$ by the projection $p: M \to K \backslash M$.  Cartan's homotopy
inverse to $p^*$ is constructed as follows.  Let
$$\alpha \in \Omega^1(M,\k)^K, \ \ \iota(\xi_M) \alpha = - \xi, \ \
\forall \xi \in \k $$
be a connection $1$-form on $M$ and
$$ \curv(\alpha) \in \Omega^2(K \backslash M,M(\k)), 
\ \ p^* \curv(\alpha) = \d \alpha + \hh [\alpha,\alpha] $$
denote its curvature.  Let
$$\pi_{\alpha} : \ \Omega(M) \to p^* 
\Omega(K \backslash M) $$
be the horizontal projection defined by $\alpha$.  The map
\begin{equation} \label{cartan}
 \Omega_K(M) \to \Omega(M)^K,
\ \ \ \eta \otimes h
\mapsto 
\left( (\pi_{\alpha} \eta) \otimes h \left( \frac{p^* \curv(\alpha)}{2 \pi i} 
\right) \right) \end{equation}
has image contained in the space of basic forms and descends to a map
$ \Omega_K(M) \to \Omega(K \backslash M)$ which is a
homotopy inverse to $p^*$.

Suppose that $M$ is compact and oriented.  Integration over $M$
vanishes on equivariant exact forms and defines a push-forward
$$ I_{M,K}: \ H_K(M) \to S(\k^*)^K .$$
We also denote by $I_{M,K}(\eta)$ the push-forward of the cohomology
class of an equivariant form $\eta$.

For any $K$-equivariant real oriented vector bundle $E$ of even
dimension $2n$, let $\Eul(E) \in H^{2n}_K(M)$ denote its equivariant
Euler class, defined as follows.  Equip $E$ with a Euclidean metric,
let $\F(E)$ denote the orthogonal frame bundle of $E$, and let $\alpha
\in \Omega^1(\F(E),\so(2n))^K$ be a $K$-invariant connection $1$-form
for $E$.  For any $\zeta \in\k$, the pairing $\alpha(\zeta_{\F(E)})$
is $K$-invariant and descends to a map
$$ \phi: \ M \to \Hom(\k,\so(E))  .$$  
The form $\curv_\k(E) \in \Omega^2_K(M,\so(E))$ defined by
$$ \curv_\k(E) := \curv(E) + 2 \pi i \phi $$
is the {\em equivariant curvature} of $E$.  The Euler class of $E$ is
$$ \Eul(E) :=  \Pf \left( \frac{ \curv_\k(E)}{2\pi} \right) \in
\Omega^{2n}_K(M)$$
where $\Pf$ is the Pfaffian, and the right-hand side denotes the
Chern-Weil characteristic form.

Let $K$ be a torus.  If $E$ is a complex $K$-representation then $E$
splits into a sum of weight spaces $E_\mu$ for $\mu \in \k^*$, so that
$ \exp(\xi) v = e^{ 2\pi i \mu(\xi)}v $ for $v \in E_\mu$ and $\xi \in
\k$.  If $E$ is a real even-dimensional representation of $K$, then
$E$ admits an invariant complex structure and the weights
$\mu_1,\ldots,\mu_n$ are independent of the choice of complex
structure up to sign.  If $E$ is oriented then the product of the
complex weights is determined by the orientation on $E$ and
$$ (\Eul(E))(\xi) = \prod_{j=1}^n - 2\pi i \mu_j(\xi).$$

We will also need cohomology with smooth and distributional
coefficients.  An equivariant form with smooth coefficients is a
smooth equivariant map $\k \to \Omega(M)$.  The equivariant
differential extends to equivariant forms with smooth coefficients and
its cohomology is the {\em equivariant cohomology of $M$ with smooth
coefficients}.  Let $C_0^\infty(\k^*)$ denote the space of compactly
supported smooth functions on $\k^*$, and $D'(\k^*)$ the space of
distributions on $\k^*$, that is, the space of continuous linear forms
on $C_0^\infty(\k^*)$.  Let $\S(\k^*)$ denote the space of {\em
Schwartz functions} on $\k^*$, the space of smooth functions $f$ such
that for any polynomial differential operator $P$, the function $Pf$
is bounded.  Its dual $\S'(\k^*)$ is the space of {\em tempered
distributions} on $\k^*$.  The inclusion $C_0^\infty(\k^*) \subset
\S(\k^*)$ dualizes to an injection $S'(\k^*) \to \D'(\k^*)$.  The
symmetric algebra $S(\k^*)$ embeds in $\S'(\k^*)$ via Fourier
transform as the space of distributions supported at the identity.  An
equivariant differential form with distributional coefficients is an
equivariant continuous linear map from $C_0^\infty(\k^*)$ to $\Omega
(M)$.  Let $\mathcal{C}_K(M)$ denote the complex of such forms; the
equivariant differential extends to $\mathcal{C}_K(M)$ and its
cohomology $\H_K(M)$ is the equivariant cohomology with distributional
coefficients.  The basic results on equivariant cohomology with
distributional coefficients are discussed in detail in Kumar-Vergne
\cite{ku:eq}.  If $M$ is compact and connected, then the map $I_{M,K}$
extends to $ I_{M,K}: \ \H_K (M) \to \D'(\k^*)^K .$

In order to state the localization formula, we need to discuss
inversion of the Euler class.  Suppose that $E$ is an oriented real
vector bundle of even dimension, and that a circle subgroup
$U(1)_\zeta \subset G$ generated by $\zeta \in \k$ acts trivially on
$M$ fixing only the zero section in $E$.  According to Atiyah-Bott
\cite{at:mo}, the Euler class is invertible after suitably modifying
the coefficient ring.  For distributional coefficients, the
construction is carried out in Paradan \cite[Section 4]{pa:loc}.  A
definition equivalent to Paradan's goes as follows.  Since $\zeta$
acts with non-zero weights $\Pf(\phi(m,\xi))$ is a {\em hyperbolic
polynomial}, that is,
$$ \Pf(\phi(m, \xi + i \tau \zeta)) \neq 0, \ \ \forall \xi \in \k,
\tau < 0, m \in M. $$
By a standard result in distribution theory \cite[Theorem
12.5.1]{ho:an2} $\Pf(\phi(m,\cdot))^{-k}$ has a unique distributional
extension with support on $(\zeta,\cdot) \ge 0$ for any $k> 0$.  Let $
\Pf(\curv_\k(E)/2\pi)_+$ denote the terms containing forms on $M$ of
positive degree, so that
$$ \Pf(\curv_\k(E)/2\pi) = \Pf(i \phi) + \Pf(\curv_\k(E)/2\pi)_+ .$$
Define
$$ \Eul(E)^{-1}_\zeta := \Pf(i \phi)^{-1} \left( 1 +
 \frac{\Pf(\curv_\k(E)/2\pi)_+}{\Pf(i \phi)} \right)^{-1} $$
interpreted via its power series expansion, which is finite since
$\Pf(\curv_\k(E))_+$ is nilpotent.

For any $\zeta \in \k$, let $M^\zeta$ denote the fixed point set of
the one-parameter subgroup $U(1)_\zeta$ generated by $\zeta$.  Fix
orientations on $M^\zeta$ and $T_{M^\zeta} M$ which induce the given
orientation on $TM | M^\zeta$.  If $E$ is a $K$-equivariant vector
bundle and $K' \subset K$ is a subgroup stabilizing a submanifold $M'
\subset M$ then $\Res^{M,K}_{M',K'} E$ denotes the restriction of $E$
to a $K'$-equivariant bundle on $M'$.  Similarly, if $\eta$ is a
$K$-equivariant cohomology form or class, we denote by
$\Res^{M,K}_{M',K'} \eta$ the restriction of $\eta$ to a
$K'$-equivariant form or class on $M'$.
\begin{theorem} (Localization for one-parameter subgroups) \label{Hloc}
For any compact oriented $K$-manifold $M$, $\eta$ an equivariant
cohomology class with smooth coefficients, and $\zeta \in \k$,
$$I_{M,K_\zeta}(\Res^{K}_{K_\zeta} \eta) = I_{M^\zeta,K_\zeta}(
\Res^{M,K}_{M^\zeta,K_\zeta} \eta \wedge \Eul(T_{M^\zeta}
M)^{-1}_\zeta).$$ \end{theorem}
For smooth values of $\Eul^{-1}_\zeta(T_{M^\zeta} M)$, Theorem
\ref{Hloc} is proved in Atiyah-Bott \cite{at:mom} and Berline-Vergne
\cite{be:ze}.  For the distributional version, see
Guillemin-Lerman-Sternberg \cite{gu:on}, Canas-Guillemin \cite{ca:ko},
and Paradan \cite[Section 5]{pa:loc}.

\section{Hamiltonian $K$-manifolds}
\label{Hamsec}

Let $T \subset K$ be a maximal torus, $\t$ its Lie algebra, and $\t^*$
its dual.  Let $\t^* \to \k^*$ be the injection whose image is the
fixed point set of $T$.  Choose a closed positive Weyl chamber $\t_+$.
Using an invariant metric on $\k$ to identify $\t^*$ with $\t$, let
$\t_+^*$ denote the image of $\t_+$ in $\t^*$; this is independent of
the choice of metric.  If $N$ is a right $K$-manifold, then by $N
\times_K M$ we mean the quotient of $N \times M$ by the $K$-action $k
(n,m) = (nk^{-1},km)$.

\subsection{Basic definitions and results} 
A {\em Hamiltonian $K$-manifold} consists of a smooth $K$-manifold
$M$, a symplectic form $\omega$ and an equivariant {\em moment map}
$\Phi: \ M \to \k^*$ satisfying
\begin{equation} \label{mommap}
 \iota(\xi_M) \omega = - d (\Phi, \xi), \ \forall \xi \in \k
 .\end{equation}
If $\omega$ is closed but degenerate the data are called a {\em
degenerate Hamiltonian $K$-manifold}.  We denote by $K_M \subset K$
the {\em principal stabilizer}, that is, the stabilizer of a generic
element in $\Phinv(\t_+^*)$.  The {\em principal orbit-type stratum}
for $K$ resp. $\k$ is the set of points $m \in M$ with $K_m$ conjugate
to $K_M$ resp. $\k_m$ conjugate to $\k_M$.  For references on the
following, see \cite{ki:con},\cite{le:co}.

\begin{theorem}  Let $M$ be a connected Hamiltonian 
$K$-manifold with proper moment map.
\begin{enumerate}
\item (Kirwan Convexity) The intersection $\Delta(M) := \Phi(M) \cap
\t_+^*$ is a convex polyhedron called the {\em moment polyhedron} of
$M$.
\item (Principal cross-section) The open face $\sigma(M)$ of $\t_+^*$
containing $\Delta(M)$ in its closure is the {\em principal face} for
$M$.  The inverse image $\Phinv(K\sigma(M))$ is an open subset of $M$
whose complement is codimension at least two, and the map
%
$ K \times_{K_\sigma} \Phinv(\sigma(M)) \to M, \ \ [k,m] \mapsto km  $
%
is a diffeomorphism onto its image.
\end{enumerate}
\end{theorem}

\noindent The {\em rank} of $M$ is the dimension of $\Delta(M)$.  $M$ is {\em
maximal rank} if and only if $\k_M$ is trivial.  

The moment map condition \eqref{mommap} is equivalent to the condition
that the {\em equivariant symplectic form}
$$ \om_K(\xi) = \om + 2 \pi i (\Phi,\xi) \in \Omega_K(M).$$
is equivariantly closed.  The {\em equivariant Liouville form} is
$$ \L := \exp(\om_K) = \exp(\om) \exp( 2\pi i (\Phi,\xi)). $$

Let $M$ be a Hamiltonian $K$-manifold with proper moment map.  The
{\em Duistermaat-Heckman measure} $\mu_{M,K}$ is the push-forward of
the measure defined by the top degree component of $\exp(\omega)$
under $\Phi$,
$$ \mu_{M,K} := \Phi_*(\exp(\omega)) \in \D'(\k^*)^K.  $$
More generally, suppose that $\eta \in \Omega_K(M)$ is closed.  If $M$
is compact, we define the {\em twisted Duistermaat-Heckman
distribution} as the Fourier transform of the push-forward of $\L
\wedge \eta$:
$$ \mu_{M,K}(\eta) = \F_\k ( I_{M,K}(\L \wedge \eta)).  $$
If $M$ is a Hamiltonian $K$-manifold $M$ with proper moment map
$\Phi$, suppose that $\xi_1,\ldots, \xi_{\dim(\k)}$ are coordinates on
$\k$ and that $\eta = \sum_I \eta_I \xi_I $ where $I$ ranges over
multisets with elements $1,\ldots,\dim(\k)$ and $ \xi_I := \prod_{i
\in I} \xi_i .$ Define
\begin{equation} \label{terms}
\mu_{M,K}(\eta) := \sum_I \partial_I \Phi_*(\exp(\omega) \wedge
\eta_I) \in \D'(\k^*)^K\end{equation}
where $\partial_I = \F_\k(\xi_I)$, the Fourier transform of $\xi_I$.
$\mu_{M,K}(\eta)$ is supported on the image of $\Phi$ and depends only
on the cohomology classes of $\om_K$ and $\eta$.  

Later we will need a variation of this construction when $M$ is a
compact Hamiltonian $K$-manifold with boundary.  Let $\eta \in
\Omega_K(M)$ be closed and $\alpha \in \Omega^1_K(M) = \Omega^1(M)^K$
an invariant one-form.  Consider the family of (possibly degenerate)
symplectic forms and moment maps
$$ \omega_s = \omega + s \d \alpha, \ \ \ (\Phi_s(m),\xi)
= (\Phi(m),\xi) + s \alpha(\xi_M) .$$
Let $\mu_{M,K,s}(\eta)$ denote the corresponding twisted
Duistermaat-Heckman distribution.

\begin{proposition} \label{bound} Let $h \in C_0^\infty(\k^*)^K$ be such that 
$\on{supp}(h) \cap \Phi_s(\partial M)$ is empty for $s \in [0,1]$.
$(\mu_{M,K,s}(\eta),h)$ is independent of $s \in [0,1]$.
\end{proposition}
\noindent This follows from the same argument as in the case without
boundary, since the relevant integrals are of forms supported on the
interior.

\subsection{Coadjoint orbits} 

The following material is mostly covered in Berline-Getzler-Vergne
\cite[Section 7.5]{be:he}.  We parametrize coadjoint orbits by their
intersection with the positive chamber:
$$ \t_+^* \to K \backslash \k^*, \ \ \ \lambda \mapsto K
 \lambda.
$$ 
A symplectic form on $K \lambda$ is defined by the
Kirillov-Kostant-Souriau formula
\begin{equation} \label{kks}
 \omega_m(\xi_M(m),\zeta_M(m)) = (m, [\xi,\zeta]).
\end{equation}
%
%
The action of $K$ on $K \lambda$ is Hamiltonian with moment map given
by the inclusion into $\k^*$.  The weights on the tangent space at a
$T$-fixed point $w \lambda$ are the roots $\alpha$ with
$(\alpha,w\lambda) > 0$.  By localization \eqref{Hloc}
\begin{equation} \label{coadj} 
(I_{K \cdot \lambda,T}(\L))(\zeta) = \sum_{[w] \in W/W_\sigma} \frac{
\exp( 2\pi i (w\lambda, \zeta ))}{ \prod_{(\alpha,w\lambda) < 0} 2 \pi i
(\alpha, \zeta ) } \end{equation}
Let $\rho$ denote the half-sum of positive roots of $\k$.  The
symplectic volume of $K \cdot \lambda$ is the equal to
\begin{eqnarray} \label{coadjvol} 
 \Vol(K \cdot \lambda) &=& (I_{K \cdot \lambda,T}(\L)) (0) \\ &=&
\lim_{t \to 0} (I_{K \cdot \lambda,T}(\L)) ( t\rho) \\ &=&\frac{
\prod_{(\alpha,\lambda) < 0} (\alpha,\lambda) }{
\prod_{(\alpha,\lambda) < 0} (\alpha, \rho )} .\end{eqnarray}
A computation at the tangent space at $\lambda$ shows that the
symplectic and Riemannian volumes are related by
%
%
%
%
%
%
\begin{equation} \label{VolGT}
 \Vol(K \cdot \lambda) = (2 \pi)^{\dim(K/K_\lambda)/2} \prod_{(\alpha,\lambda) > 0}  (\alpha,\lambda) \Vol(K/K_\lambda)  \end{equation}
which implies
\begin{equation} \label{VolGT2}
 \Vol(K/K_\lambda)^{-1} = (2 \pi)^{\dim(K/K_\lambda)/2} 
\prod_{(\alpha,\lambda) > 0}  (\alpha,\rho) .\end{equation}
Similarly, let $\RVol(K \cdot \lambda)$ denote the volume of $K \cdot
\lambda$ with respect to the Riemannian metric induced by the
embedding $K \cdot \lambda \to \k^*$. We have
\begin{equation}\label{Rvol}
 \RVol(K \cdot \lambda) = (2 \pi)^{\dim(K/K_\lambda)} 
 \prod_{(\alpha,\lambda) > 0} (\alpha,\lambda)^2 \Vol(K/K_\lambda) 
\end{equation}
If $K_\lambda = T$ then 
$\RVol(K \cdot \lambda) = \Pi(\lambda)^2 \Vol(K/K_\lambda)$
where
$$\Pi(\xi) = \prod_{\alpha > 0} 2\pi (\alpha,\xi) =
i^{-\dim(K/T)/2} \Eul(\k/\t).  $$
\subsection{Symplectic quotients}

The {\em symplectic quotient} of $M$ at $\lambda \in \k^*$ is
$$ M_{(\lambda)} := K_\lambda \backslash \Phinv(\lambda) .$$
If $\Phinv(\lambda)$ is contained in the principal orbit-type stratum
for $K$ (resp. $\k)$ then $M_{(\lambda)}$ has the structure of a
symplectic manifold (resp. orbifold), with symplectic form
$\omega_{(\lambda)}$ the unique two-form that pulls back to the
restriction of $\omega$ to $\Phinv(\lambda)$.  For $\lambda$ such that
$\Phinv(\lambda)$ is contained in the principal orbit-type stratum, we
denote by $\L_{(\lambda)} = \exp(\omega_{(\lambda)})$ the Liouville
form on $M_{(\lambda)}$.

The relation between the cohomology of $M$ and the cohomology of the
quotient was studied by Kirwan \cite{ki:coh}.  Let $\lambda$ be a
regular value of $\Phi$ and $\kappa_\lambda$ the composition of
restriction $H_K(M) \to H_K(\Phinv(\lambda))$ with the isomorphism
$ H_K(\Phinv(\lambda)) \to H(M_{(\lambda)}) .$
$\kappa_\lambda$ extends to cohomology with smooth coefficients.  By
\cite{ki:coh}, if $\lambda$ is central and $M$ is compact then
$\kappa_\lambda$ is surjective.

The cohomological pairings on the symplectic quotient are encoded in
the twisted Duistermaat-Heckman distributions.  Let
$\mu_{\k^*},\mu_{\t^*}$ denote the Lebesgue measures induced by the
metrics on $\k^*,\t^*$ respectively.
\begin{proposition} \label{red}
Let $M$ be a (possibly degenerate) Hamiltonian $K$-manifold with
proper moment map $\Phi$.  Let $\eta \in \Omega_K(M)$ be a closed
equivariant form.  Let $U \subset \k^*$ be a subset such that $K$ acts
locally freely on $\Phinv(U)$, and $U^{\reg} \subset U$ the set of
regular values of $\Phi$ in $U$.
\begin{enumerate}
\item 
$ \mu_{M,K}(\eta)|_U = {\Vol(K/K_M)}{\Vol(K \cdot \lambda)}^{-1}
I_{M_{(\lambda)}} ( \L_{(\lambda)} \wedge \kappa_\lambda(\eta))
\mu_{\k^*} |_U$.
\item 
$I_{M_{(\lambda)}} ( \L_{(\lambda)} \wedge \kappa_\lambda(\eta))$ is
equal to a polynomial in $\lambda$ on any subset of $U^{\reg}$ on
which $K_\lambda$ is constant.
\item 
Let $\lambda \in \t_+^*$ be a regular value of $\Phi$.  For $\nu \in
\t^*_{\on{reg}} \cap \t_+^*$,
$$ \lim_{\nu \to \lambda} (\# W_\lambda)^{-1} \I_{M_{(\nu)}} (
\L_{(\nu)} \wedge \kappa_\nu(\eta \wedge \Eul((\k_\lambda/\t)^*))) =
I_{M_{(\lambda)}} ( \L_{(\lambda)} \wedge \kappa_\lambda(\eta)) .$$
\end{enumerate}
\end{proposition}
\begin{proof}  In the non-degenerate case, (a),(b) are basic results of 
Duistermaat-Heckman, see \cite{du:eq},\cite{je:lo1} \label{Hred}.  The
degenerate cases are similar, using the machinery of equivariant
cohomology instead of local models, as in Atiyah-Bott \cite{at:mom}:
If $\lambda_1,\lambda_2 \in U^{\on{reg}}$ and $K_{\lambda_1} =
K_{\lambda_2}$ then let $A \subset \k^*$ be a one-manifold (arc) with
boundary $\partial A = \{ \lambda_1,\lambda_2 \}$, such that $A$ is
contained in $U$ and $K_\lambda - K_{\lambda_1}$ for all $\lambda \in
A$.  Suppose that the inverse image $\Phinv(A)$ is smooth.  Since $K$
acts locally freely, the quotient
$$ M_A = K_{\lambda_1} \backslash \Phinv(A) $$
is a smooth orbifold.  By \eqref{cartan}, the equivariant symplectic
form $\ti{\om}$ as well as $\eta$ descend to closed forms on $M_A$.
We have
\begin{equation} \label{Lk}
 \L_{\lambda_i} = \kappa_{\lambda_1}( \exp(\ti{\om} - 2 \pi i
 (\lambda_i,\xi))), \ \ i = 1,2 \end{equation}
and therefore by Stokes theorem applied to $M_A$
\begin{equation} \label{Xform}
I_{M_{(\lambda_2)}}(\L_{\lambda_2} \wedge \kappa_{\lambda_2} (\eta)) =
I_{M_{(\lambda_1)}}(\L_{\lambda_1} \wedge \kappa_{\lambda_2} (\eta
\wedge \exp( 2\pi i (\lambda_1 - \lambda_2,\xi))) \end{equation}
which is polynomial in $\lambda_2$. Since $\Phinv(A)$ is smooth for
generic $A$, this shows (b).  (c) $M_{(\nu)}$ is a $K_\lambda/T$-fiber
bundle over $M_{(\lambda)}$.  The class
$\kappa_\lambda(\Eul((\k_\lambda/\t)^*))$ restricts to the Euler class
of the tangent bundle on the fiber $K_\lambda/T$, which has Euler
characteristic $\# (K_\lambda/T)^T = \# W_\lambda$.  The result
follows from fiber integration.
\end{proof}  

We will need a more general ``reduction in stages'' version of this
result.  Let $K_1 \subset K$ be a normal subgroup, and $K_2 = K/K_1$
the quotient.  Let $\Phi = (\Phi_1,\Phi_2)$ be the decomposition of
$\Phi$ according to an invariant splitting $\k \cong \k_1 \oplus
\k_2$.
\begin{proposition} \label{stages} Let $M$ be a (possibly degenerate)
Hamiltonian $K$-manifold with proper moment map $\Phi$.  Let $\eta \in
\Omega_K(M)$ be closed and $U_1 \subset \k_1^*, U_2 \subset \k_2^*$
open subsets such that $K_1$ acts locally freely on $\Phi^{-1}(U_1
\times U_2)$.  Then
$$ \mu_{M,K}(\eta) |_{U_1 \times U_2} =
\frac{\Vol(K_1/K_{1,M})}{\Vol(K_1 \cdot \lambda_1)} \int_{U_1^{\reg}}
( \mu_{M_{(\lambda_1)},K_2}(\kappa_{\lambda_1}(\eta)) |_{U_2} \otimes
\delta_{\lambda_1}) \d \lambda_1 $$
where the integral is over the set $U_1^{\reg}$ of regular values of
$\Phi_1$ in $U_1$.  \end{proposition}
\noindent Here the twisted Duistermaat-Heckman distribution $
\mu_{M_{(\lambda_1)},K_2}(\eta_{(\lambda_1)})$ on $M_{(\lambda_1)}$ is
a distribution on $\k^*_2$.  Its tensor product with
$\delta_{(\lambda_1)}$ is a distribution on $\k^*$ depending on
$\lambda_1$.  

\begin{lemma}  \label{stagepol} Suppose that, in the setting of Proposition 
\ref{stages}, $\k_m \subset \k_2$ for all $m \in \Phinv(U_1 \times
U_2)$.  Then $ (\mu_{M_{(\lambda_1)},K_2}(\kappa_{\lambda_1}(\eta))
|_{U_2}$ depends polynomially on $\lambda_1$.
\end{lemma}

\begin{proof}  The assumption $\k_m \subset \k_2$ implies that the 
$K_{1,\lambda_1}$-bundle
$ \Phi_1^{-1} (\lambda_1) \to M_{(\lambda_1)} $
admits a connection one-form $\alpha$ that vanishes on the generating
vector fields for $K_2$.  (Construct the connection locally, then
patch together.)  Hence the $K_2$-equivariant form $\curv_{K_2}\alpha$
is nilpotent, so
$ \kappa_{\lambda_1}(\exp(2 \pi i (\lambda_1,\xi))) = \exp((
 \lambda_1, \curv_{K_2}\alpha)) $
is a polynomial in $\lambda_1$.  Using the $K_2$-equivariant analog of
\eqref{Lk} implies the result.
\end{proof}  

We apply this to prove a polynomiality result for non-regular values.
Let $\lambda_0,\lambda_1 \in \k^*$ be central and
$$R_{\lambda_0,\lambda_1} = \lambda_0 + \R_{\ge 0} (\lambda_1 -
\lambda_0)$$
the ray starting from $\lambda_0$ to $\lambda_1$.  Let $\k_1$ denote
the span of $\lambda_1$, and $\k_2$ the quotient $\k/\k_1$.  We say
that a distribution $\mu \in \D'(\k^*)$ is smooth, resp. polynomial
along the ray $R_{\lambda_0,\lambda_1}$ at $\lambda_2 \in
R_{\lambda_0,\lambda_1} - \{ \lambda_0 \}$ if $\mu$ is equal to
$$ \int_{\k_1^*} (\mu_2(\lambda_1) \otimes \delta(\lambda_1)) \d
\lambda_1 $$
for some distribution $\mu_2(\lambda_1) \in \D'(\k_2^*)$ depending
smoothly, resp. polynomially on $\lambda_1 \in \k_1^*$.  We say $\mu$
is smooth, resp. polynomial near $\lambda_2 \in
R_{\lambda_0,\lambda_1}$ if this holds in a neighborhood of
$\lambda_2$.

\begin{proposition} \label{ray1} Let  $M$ be a (possibly degenerate)
Hamiltonian $K$-manifold with proper moment map $\Phi$.  Let $\eta \in
\Omega_K(M)$ be closed.  Let $\lambda_0,\lambda_1 \in \k^*$ be
distinct and central.  Then $\mu_{M,K}(\eta)$ is polynomial along
$R_{\lambda_0,\lambda_1}$ near $\lambda_1$ sufficiently close to
$\lambda_0$.
\end{proposition}  

\begin{proof}   Let $m \in \Phinv(\lambda_1)$, $\xi \in \k_m$, and 
$M^\xi$ the infinitesimal fixed point set of $\xi$ containing $m$.
For $\lambda_1$ sufficiently close to $\lambda_0$, $M^\xi$ meets
$\Phinv(\lambda_0)$ and so
$$ (\lambda_1,\xi) = (\Phi(M^\xi),\xi) = (\lambda_0,\xi) .$$  
It follows that $\lambda_1 - \lambda_0 $ annihilates $\xi$, hence
$\xi$ is contained in $K_2$.  The result now follows from Proposition
\ref{stages} and Lemma \ref{stagepol}.
\end{proof}

\subsection{Induction}

First we define induction for distributions.  Let $\tau$ be any face
of $\t_+^*$.  Let $\Vol^K_{K_\tau}: \t^* \to \R$ denote the function
$$ \Vol^K_{K_\tau}(\lambda):= \frac{\Vol(K \cdot \lambda)}{\Vol(K_\tau
\cdot \lambda)}, \ \ \lambda \in \t_+^* .$$
Using \eqref{coadjvol} $\Vol^K_{K_\tau}$ has a polynomial extension to
$\t^*$ that is invariant under $W_\tau$. We denote by the same name
its extension to $\k_\tau^*$.  Define
\begin{equation} \label{Ind}
 \Ind_{K_\tau}^K : \ \D'(\k_\tau^*) \to \D'(\k^*)^K, \ \ \
(\Ind_{K_\tau}^K \mu,h) = ( \mu, \Vol^K_{K_\tau} \Res_{K_\tau}^K h). 
\end{equation}
Restriction to tempered distributions defines a map $\S'(\k_\tau^*)
\to \S'(\k^*)^K$.  The same notation will be used for the Fourier
transform $\S'(\k_\tau) \to \S'(\k)^K .$ The reader may note that
$\Ind_{K_\tau}^K$ is the ``semiclassical limit'' of holomorphic
induction of representation rings $ R(K_\tau) \to R(K) .$

Next, we define induction for Hamiltonian actions.  If $M$ is a
Hamiltonian $K_\tau$-manifold, one can define a Hamiltonian
$K$-manifold by
$$ \Ind_{K_\tau}^K M := K \times_{K_\tau} M $$
with the unique closed equivariant two-form $\Ind_{K_\tau}^K
\om_K$ restricting to $\om_K$ on $M$.  The two-form
$\Ind_{K_\tau}^K \omega$ is degenerate if and only if $\Delta(M)$ lies
in the union of open faces of $\t_+^*$ whose closure contains $\tau$.

Finally, we define induction for equivariant forms.  The inclusion
$$ M \to \Ind_{K_\tau}^K M, \ \ \ m \mapsto [1,m] $$
induces a map 
%
$ \Omega_{K}(\Ind_{K_\tau}^K M) \to \Omega_{K_\tau}(M) .$
%
A homotopy inverse is provided by the composition $\Ind_{K_\tau}^K$ of
the maps
\begin{equation} \label{cartan2}
 \Omega_{K_\tau}(M) \to \Omega_{K \times K_\tau}(K
\times M) \to \Omega_{K}(\Ind_{K_\tau}^K M) \end{equation}
where the last map is the Cartan map \eqref{cartan}.  

The following proposition shows that taking Duistermaat-Heckman
distributions commutes with induction, see Paradan
\cite[3.13]{pa:mo}. For completeness we include a proof.

\begin{proposition} \label{cross2}   Let $\tau$ be a face of the
  positive Weyl chamber, $M$ a Hamiltonian $K_\tau$-manifold with
proper moment map, and $\eta$ a closed polynomial $K_\tau$-equivariant
form on $M$.  $ \mu_{\Ind_{K_\tau}^K M,K}(\Ind_{K_\tau}^K \eta) =
\Ind_{K_\tau}^K \mu_{M,K_\tau}(\eta) .$
\end{proposition} 
\begin{proof}  If the closure of $\tau$ contains
the principal face $\sigma$ of $M$, then 
$ \Ind_{K_{\sigma}}^K = \Ind_{K_\tau}^K \Ind_{K_\sigma}^{K_\tau} .$
Therefore, it suffices to prove the proposition for $\tau = \sigma$.
Let $ \on{conn}^K_{K_\sigma} \in \Omega^1(K,\k_\sigma)$ be the
connection on $K \to K/K_\sigma$ defined using the metric on $\k$, $
\widetilde{\curv}^K_{K_\sigma} \in
\Omega^2_K(K/K_\sigma,K(\k_\sigma))$ its equivariant curvature and
$\curv^K_{K_\sigma} \in \Omega^2(K/K_\sigma(\k_\sigma)$ its ordinary
curvature.  For each $\lambda \in \sigma$ the pairing with the
curvature defines an equivariant two-form $
(\widetilde{\curv}^K_{K_\sigma},\lambda) \in \Omega^2_K(K/K_\sigma) .$
Let
$$ p_i \in S(\k_\sigma)^{K_\sigma}, i = 1,2,\ldots $$
be a basis for the invariant polynomials on $\k_\sigma$.  For each $i$
we have a characteristic form defined via the Chern-Weil homomorphism
$$ \widetilde{\curv}^K_{K_\sigma,i} = \sum_I {\curv}^K_{K_\sigma,i,I}
\xi_I \in \Omega_K(K/K_\sigma) .$$
Because $(\curv^K_{K_\sigma},\lambda)$ is the pull-back of
the Kirillov-Kostant-Souriau form \eqref{kks} under the map
$K/K_\sigma \to K \cdot \lambda$, we have for any $h \in \S(\k^*)^K$
\begin{equation} \label{triv}  \sum_I \int_{K/K_\sigma} \curv^K_{K_\sigma,i,I} \wedge
\exp(\curv^K_{K_\sigma},\lambda) (\partial_I h)(\lambda) = (\partial_i
\Vol_{K_\sigma}^K h) (\lambda) \end{equation}
where $\partial_i$ is the Fourier transform of $p_i$.  Let $\eta_i \in
\Omega(M)$ be forms such that
$\eta = \sum_i \eta_i p_i .$
$ \Ind_{K_\sigma}^K (\eta)$ is the form on $\Ind_{K_\sigma}^K M$ whose
pull-back to $K \times M$ is
$$ \sum_i \pi_2^*  \eta_i \wedge \pi_1^* \phi^* (
\curv^K_{K_\sigma})_i .$$
Let $\beta \in \Omega(K \times M)$ be a form which integrates to $1$
on the orbits of $K_\sigma$ on $K \times M$.  Using \eqref{triv} we
have (omitting pull-backs which confuse the notation)
\begin{eqnarray*} 
(\mu_{\Ind_{K_\sigma}^K M,K}(\Ind_{K_\sigma}^K(\eta)),h) &=&
\sum_{i,I} \int_{K \times M} (\Ind_{K_\sigma}^K \Phi)^* (\partial_I h)
\exp(\om) \wedge\eta_i \wedge \\ && \curv^K_{K_\sigma,i,I} \wedge
\beta \wedge \exp(\curv^K_{K_\sigma},\Ind_{K_\sigma}^K \Phi) \\
&=& \sum_{i}\int_{M}  \Phi^* ( \partial_i \Vol_{K_\sigma}^K \Res_{K_\sigma}^K
h)  \exp(\om) \wedge \eta_i  
\\
&=& (\mu_{M,K_\sigma}(\eta), \Vol_{K_\sigma}^K  \Res_{K_\sigma}^K h) \\
&=& (\Ind_{K_\sigma}^K \mu_{M,K_\sigma}(\eta),h).
\end{eqnarray*}
as claimed.
\end{proof}

For $\lambda_0,\lambda_1 \in \k^*$ distinct, we say that a
distribution $\mu$ on $\k^*$ is smooth, resp. polynomial along the ray
$R_{\lambda_0,\lambda_1}$ at $\lambda \in R_{\lambda_0,\lambda_1} - \{
\lambda_0 \}$ if $\mu$ is the induction of a distribution on
$\k_\lambda^*$ that is polynomial along $R_{\lambda_0,\lambda_1}$ near
$\lambda$.  By Proposition \ref{cross2}, Proposition \ref{ray1} holds
without the assumption that $\lambda_1$ is central.

\subsection{Comparison of abelian and non-abelian Duistermaat-Heckman 
distributions}

For the sake of computing examples it will be helpful to have the
formula that compares the abelian and non-abelian Duistermaat-Heckman
measures.  The following result of Harish-Chandra compares Fourier
transforms over $\k$ and $\t$:
\begin{lemma} \label{hc} For any $\ h \in \S(\k)^K$, 
%
$ \Pi \Res^{\k^*}_{\t^*}\F_\k(h) = i^{\dim(K/T)/2} \F_\t(\Pi
 \Res^{\k}_{\t} h). \ \
$
%
\end{lemma}

\begin{proof} Let $\lambda \in \t^*$. 
\begin{eqnarray*}
(\F_\k(h))(\lambda) &=& 
(2\pi )^{-\dim(\k/2)} \int_{\k} e^{2 \pi i (\lambda,\xi)} h(\xi)
\d \xi \\ 
&=& 
(2\pi )^{-\dim(\k/2)} \int_{\t_+ \times K/T} e^{2 \pi i (\lambda,k\xi)} h(k
\cdot \xi) \exp(\omega_\xi) \frac{\RVol_T^K(\xi)}{\Vol_T^K(\xi)} \d \xi \\
&=& (2 \pi )^{-\dim(\t)/2} 
\int_{\t_+} \sum_{w \in W} (-1)^{l(w)} \frac{h(\xi) e^{2 \pi i
(\lambda,w\xi)} \Pi(\xi) 
}{\Eul(\k/\t)(\lambda)} \d \xi \\ 
&=& (2 \pi)^{-\dim(\t)/2}  i^{\dim(\k/\t)/2}  \int_{\t} \frac{h(\xi)
  e^{2 \pi i (\lambda,\xi)} \Pi(\xi)}{\Pi(\lambda)} \d \xi \\ 
&=& i^{\dim(\k/\t)/2} (\Pi^{-1} \F_\t ( \Res^K_T h \cdot \Pi)
) ({\lambda}) .\end{eqnarray*}
\end{proof}

As I learned from P.-E. Paradan, Harish-Chandra's result implies the
following relation between Duistermaat-Heckman measures.  Note that
the Euler class $\Eul(\k/\t)$ considered as a distribution on $\t^*$
is the product of partial derivatives in the direction of the negative
roots of $\k$.

\begin{theorem} \label{Sin} Let $M$ be a compact Hamiltonian
 $K$-manifold and $\eta \in \Omega_K(M)$ closed.
\label{TtoG} $ \mu_{M,K}(\eta) = (\#W)^{-1} \Ind_T^K \Eul(\k/\t)
\mu_{M,T}(\Res^K_T \eta) .$
\end{theorem} 

\begin{proof}  Since $ \Res_T^K
\F_\k^{-1} \mu_{M,K}(\eta) = \F_\t^{-1} \mu_{M,T}(\eta) $ we have
using Lemma \ref{hc}
\begin{eqnarray*} 
 (\mu_{M,K}(\eta), \F_\k(h))
& = & (\F_{\k}^{-1} ( I_{M,K}(\L \wedge \eta)), h) \\
& = & (\# W)^{-1} (\RVol_T^K \F_{\t}^{-1} (I_{M,T}(\L \wedge \eta)),
\Res_T^K h) \\
& =& (\# W)^{-1} (I_{M,T}(\L \wedge \eta), \F_{\t}( \RVol_T^K \Res_T^K
h )) \\
&= & (\# W)^{-1} ( I_{M,T}(\L \wedge \eta), \F_\t( \Pi^2 \Vol(K/T)
\Res_T^K h) ) \\
& = & (\# W)^{-1} ( \mu_{M,T}(\eta), \Eul(\k/\t) \Pi 
\Vol(K/T) \Res_T^K \F_\k ( h) ) \\
& = & (\# W)^{-1} ( \mu_{M,T}(\eta), \Eul(\k/\t) 
\Vol_T^K \Res_T^K \F_\k ( h) ).
\end{eqnarray*}
\end{proof}

This formula has as a corollary a result of S. Martin \cite{mar:sy},
which compares cohomological pairings on the abelian and non-abelian
quotients.  We denote by $M_{T,(\lambda)}$ the symplectic quotients
for the action of $T$:
$$ M_{T,(\lambda)} = (\Res_T^K \Phi)^{-1}(\lambda)/T .$$
\begin{proposition} If $\lambda$ is a regular value of $\Phi$ 
and $\Res_T^K \Phi$ then
$$ I_{M_{(\lambda)}}(\kappa_\lambda(\eta)) = (\# W_\lambda)^{-1}
I_{M_{(\lambda),T}}(\kappa_{\lambda}(\Res_T^K \eta \wedge \Eul(\k/\t)
\wedge \Eul((\k/\k_\lambda)^*))).$$
\end{proposition}  

\begin{proof}  From \eqref{Sin} and \ref{red} (a) we have for 
generic $\lambda$
$$  I_{M_{(\lambda)}}(\kappa_\lambda(\eta))
= I_{M_{(\lambda),T}}(\kappa_{\lambda}(\Res_T^K \eta \wedge \Eul(\k/\t)
)) .$$
The result for arbitrary $\lambda$ follows from \ref{red} part
(c).
\end{proof}  

\subsection{Symplectic vector bundles} 

Let $M$ be a compact Hamiltonian $K$-manifold and $\pi: \ E \to M$ a
$K$-equivariant symplectic vector bundle, that is, a vector bundle
with structure group $Sp(2n,\R)$.  We recall from \cite{gu:symu} that
the total space of $E$ can be given the structure of closed two-form,
equal to $\omega$ on the zero section and non-degenerate in a
neighborhood of it: Let $\Fr(E)$ denote the frame bundle of $E$ and
$\omega_F$ the symplectic form on the fiber $F := \R^{2n}$.  The
action of $\Sp(2n,\R)$ on $F$ is Hamiltonian; we denote by $\phi: \, F
\to \lie{sp}(2n,\R)^*$ the moment map. Let $\alpha \in
\Omega^1(\Fr(E),\lie{sp}(2n,\R))$ be a connection one-form. The two-form
$$\pi^*\om + d(\alpha,\phi) + \om_F \in \Omega^2(\Fr(E) \times F)$$ 
(pullbacks from factors are omitted from the notation) is basic and
descends to a closed form $\om_E$ on $E \cong \Fr(E)
\times_{Sp(2n,\R)} F$ with the required properties.  By the symplectic
embedding theorem, $\om_E$ is the unique form with these properties up
to symplectomorphism on a neighborhood of the zero section. The
construction also works equivariantly: If $M$ is a Hamiltonian
$K$-manifold and $E$ a $K$-equivariant symplectic vector bundle let
$\Phi_F$ denote the moment map for the $K$-action on the fiber $F$.
The map
$$\xi \mapsto \pi^* (\Phi,\xi) + (\alpha(\xi_{\Fr(E)}),\phi) $$ 
is $\Sp(2n,\R)$-invariant and descends to a moment map $\Phi_E: E \to
\k^*$.

Let $U(1)_\zeta$ be the one-parameter subgroup generated by a central
element $\zeta \in \k$.  Suppose $U(1)_\zeta$ acts on $E$ fixing only
the zero section with positive weights.  In this case the moment map
for the action of $U(1)_\zeta$ on the fiber $F$ is a positive-definite
quadratic form; it follows that $(\Phi_E,\zeta)$ is proper, so
$\Phi_E$ is proper as well.  For any closed form $\eta \in
\Omega_K(M)$, localization applied to the total space of $E$ gives
\begin{equation}  
\mu_{E,K}(\pi^* \eta) =
\mu_{M,K}(\eta \wedge \Eul^{-1}_\zeta(E)).
\end{equation} 
Non-compactness of $E$ can be remedied as in \cite{pr:dh} or
\cite{pa:mo}.

Suppose that $U(1)_\zeta$ acts on $E$ with both positive and negative
(but not zero) weights.  Let $E = E_- \oplus E_+$ be the decomposition
into positive and negative weight bundles.  Let $\pi: \ E' \to M$ be
the symplectic vector bundle obtained from $E$ by reversing the
symplectic structure on the subbundle $E_- \subset E$ on which
$U(1)_\zeta$ acts with negative weights, so that the orientation of 
$E'$ is $(-1)^{\dim(E_-)}$ times the orientation on $E$.  Since
$$ \Eul(E)_\zeta^{-1} = (-1)^{\dim(E^-)} \Eul(E')_\zeta^{-1} $$
we have 
\begin{equation} \label{vb}
 \mu_{M,K}(\eta \wedge \Eul(E)_\zeta^{-1}) = (-1)^{\dim(E^-)}
\mu_{E',K}( \pi^* \eta).\end{equation}
By Proposition \ref{ray1} applied to $E'$, for any $\lambda_0 \in
\k^*$ and non-zero $\lambda_1 \in \k_{\lambda_0}^*$, the distribution
$\mu_{M,K}(\eta \wedge \Eul(E)_\zeta^{-1})$ is polynomial along
$R_{\lambda_0,\lambda_1}$ near $\lambda_0$.

\subsection{Local Normal Form}

Let $(M,\omega,\Phi)$ be a Hamiltonian $K$-manifold and $y \in
\Phinv(\xi)$.  Let $\k_{y}^\circ \subset \k^*$ denote the annihilator
of $\k_{y}$.  Let $N$ denote the orthogonal complement of the tangent
space to the orbit $T_{y} (K \cdot y)$.  Let $S = N/(N^\omega \cap N)$
denote the quotient of $N$ by the kernel of the symplectic form
restricted to $N$, called the {\em symplectic slice} at $y$.  Let
$\omega_S$ denote the two-form on $S$, and $\Phi_S$ the quadratic
moment map for the action of $K_y$.  The $K$-manifold
\begin{equation} \label{model}
 M_0 = K \times_{K_{y}} (S \oplus (\k_\xi^* \cap \k_{y}^\circ)) .\end{equation}
has a Hamiltonian $K$-structure with moment map
\begin{equation} \label{momentf} \Phi_0: \ [k,s,\nu] \mapsto k \cdot (
  \Phi_S(s) + \nu + \xi) .\end{equation}
The following is proved in Marle \cite{ma:vo} and 
Guillemin-Sternberg \cite{gu:no}:
\begin{theorem} 
There exists an equivariant symplectomorphism $\psi$ of a neighborhood
$U$ of $y$ in $M$ with a neighborhood $U_0$ of $[1,0,0]$ in $M_0$.
\end{theorem} 
\noindent 

\section{The Kirwan-Ness stratification}
\label{Kirsec}

Let $(M,\omega)$ be a Hamiltonian $K$-manifold with proper moment map
$\Phi$ and $f$ one-half the norm-square of the moment map,
$$f: \ M \to \R, \ f(m) = \hh (\Phi(m),\Phi(m)).$$ 
In general, $f$ is not a Morse-Bott function.  The critical set of $f$
consists of points $m$ fixed by the vector field generated by
$\Phi(m)$:
\begin{equation} \label{critf}
 \crit(f) = \{ m \in M, \ \ (\Phi(m)_M)(m) = 0 \}. \end{equation}
Hence $\Phinv(0)$ is a component of $\crit(f)$.  For any connected
component $C \subset \crit(f)$ the intersection $\Phi(C) \cap \t_+$
consists of a single point $\xi$.  (See \cite[3.15]{ki:coh} for the
case $M$ compact; the case $\Phi$ is proper is similar.)  Define
$$ \Xi(M) = \{ \xi(C), \ \ C \subset \crit(f) \} .$$
For any $\xi \in \Xi(M)$ let 
$$ C_\xi = \{ m \in M, \ \ \Phi(Km) \cap \t_+ = \xi \} $$
which may be a finite union of connected components.  Choose a
$K$-invariant almost complex structure on $M$, and consider the
corresponding $K$-invariant Riemannian metric.  For any $m \in M$ let
$\{ m_t , \ t \in [0,\infty)\}$ denote the trajectory of $-\grad(f)$.
By Theorem \ref{norm} of the appendix, $m_t$ converges to a
critical point of $f$ as $t \to \infty$.  For any $\xi \in \Xi(M)$,
let $ M_\xi$ denote the stable set of the corresponding critical
component $C_\xi$, that is, the set of $m$ with limit point in
$C_\xi$.  The Kirwan-Ness stratification is
$$ M = \bigcup_{\xi \in \Xi(M)} M_\xi .$$
For each $\xi \in \Xi(M)$ let $U(1)_\xi$ denote the one-parameter
subgroup generated by $\xi$.  Let $Z_\xi$ denote the union of
components of the fixed point set of $U(1)_\xi$ meeting $C_\xi \cap
\Phinv(\t_+)$, $Y_\xi$ the set of points in $M$ which flow to $Z_\xi$
under $-\grad(\Phi,\xi)$, and $ \varphi_\xi: \ Y_\xi \to Z_\xi$ the
map given by the limit of the flow.  By \eqref{critf},
$$Z_\xi \cap \Phinv(\xi) = C_\xi \cap \Phinv(\xi) .$$
Let $G$ denote the complexification of $K$, $K_\xi,G_\xi$ the
stabilizers of $\xi$ under the adjoint action of $K,G$ and $P_\xi$ the
standard parabolic corresponding to $\xi$.  Since $(\Phi,\xi)$ is a
Morse-Bott function, $Y_\xi$ is a smooth $K_\xi$-invariant
submanifold.  Let $Z_\xi^\circ$ denote the set of points in $Z_\xi$
which flow to $C_\xi$ under $-\grad(\Res^M_{Z_\xi} f)$, and
$Y_\xi^\circ = \varphi_\xi^{-1}(Z_\xi^\circ)$.  By the stable manifold
theorem (see e.g. \cite{sh:gl}) there exists a diffeomorphism
\begin{equation} \label{fiber1}
 Y_\xi^\circ = T_{Z_\xi^\circ} Y_\xi^\circ \end{equation}
where $T_{Z_\xi^\circ} Y_\xi^\circ$ is the normal bundle of
$Z_\xi^\circ$ in $Y_\xi^\circ$.  The following combines results from
Kirwan \cite[4.16,4.17]{ki:coh}, Ness \cite{ne:st}, and Heinzner-Loose
\cite{he:re}.

\begin{theorem} \label{knhl}  Let $M$ be a Hamiltonian
$K$-manifold with proper moment map, and $\xi \in \Xi(M)$.
\begin{enumerate}
\item For a suitable choice of invariant almost complex structure, the
stratum $M_\xi$ is a smooth invariant submanifold which is identical
in a neighborhood of $C_\xi$ to $KY_\xi^\circ$;
\item Suppose that $M$ is equipped with an invariant K\"ahler
structure.  For the metric defined by the structure $M_\xi$ is a
$G$-invariant K\"ahler submanifolds, $Y_\xi^\circ$ is $P_\xi$-stable
and there exist equivariant diffeomorphisms
\begin{equation} \label{fiber2}
K \times_{K_\xi} Y_\xi^\circ \to G \times_{P_\xi} Y_\xi^\circ 
\to   M_\xi.\end{equation}
\end{enumerate}
\end{theorem}

\section{Localization for the norm-square of the moment map}

Define $ \nu_\xi := T_{M_\xi} M |_{Z_\xi^\circ} \oplus T_{Z_\xi^\circ}
Y_\xi^\circ $.  In this section we will prove
\begin{theorem} \label{Hform} Let $M$ be a Hamiltonian $K$-manifold
with proper moment map and $\eta \in \Omega_K(M)$ closed.  The
restriction of $\mu_{Z_\xi,K_\xi}( \Res^{M,K}_{Z_\xi,K_\xi} \eta
\wedge \Eul(\nu_\xi)_\xi^{-1}) $ to a neighborhood of $\xi$ has a
unique extension $\mu_{Z_\xi^\circ,K_\xi}(
\Res^{M,K}_{Z_\xi^\circ,K_\xi} \eta \wedge \Eul(\nu_\xi)_\xi^{-1}) $
that is polynomial on any ray starting at $\xi$, and
\begin{equation} \label{Fform}
 \mu_{M,K} (\eta ) = \sum_{\xi \in \Xi(M)} \Ind_{K_\xi}^K
\mu_{Z_\xi^\circ,K_\xi}( \Res^{M,K}_{Z_\xi^\circ,K_\xi} \eta \wedge
\Eul(\nu_\xi)_\xi^{-1}) .\end{equation}
\end{theorem}
\noindent Properness of the moment map insures that the sum on the
right hand side of \eqref{Fform} is locally finite and so
well-defined.  See Section 7 for examples.

\subsection{Witten's deformation}
\label{deform}

Let $M$ be a Hamiltonian $K$-manifold with proper moment map $\Phi$, $
X_f \in \Vect(M)$ the Hamiltonian vector field for $f = \hh
(\Phi,\Phi)$, $g$ an invariant compatible metric on $M$, $J$ the
associated almost complex structure, and $\alpha$ the invariant
one-form
$$ \alpha(\cdot ) = g(X_f,\cdot ) = \omega(X_f, J (\cdot)).$$
We write
$$ \d_K \alpha (\xi) = d\alpha + 2\pi i (\phi,\xi), \ \ \ \ (\phi,\xi)
:= \iota(\xi_M) \alpha.$$
Note that 
\begin{equation}  \label{pos}
(\phi,\Phi) = g(X_f, X_f) \ge 0 \end{equation}
and equality holds only if $X_f = 0$. Define
$$ \omega_s := \omega + s \d \alpha, \ \ \ \ \Phi_s := \Phi + s \phi,
\ \ \ \ \ \ti{\omega_s} := \omega_s + 2\pi i \Phi_s .$$

\begin{lemma}  $\Phi_s$ is proper for all $s \ge 0$.
\end{lemma}

\begin{proof}   By \eqref{pos},
\begin{equation} \label{pos2}
\Vert \Phi_s \Vert^2 = \Vert \Phi \Vert^2 + 2s (\Phi, \phi) + s ^2
\Vert \phi \Vert^2 \ge \Vert \Phi \Vert^2.\end{equation}
Hence $\Vert \Phi_s \Vert^2 \leq C$ implies $\Vert \Phi \Vert^2 \leq
C$ which shows that $\Phi_s$ is proper. \end{proof}

Define $\mu_{M,K,s}(\eta) \in \D'(\k^*)^K$ by
\begin{equation} \label{mus}
 (\mu_{M,K,s}(\eta),h) = \sum_I \int_{M} \eta_I \wedge \exp(\omega_s)
 \partial_I h(\Phi(m) + s \phi(m)) .\end{equation}
Since the cohomology class of $\omega_s$ is independent of $s$, so is
$\mu_{M,K,s}(\eta)$.  Let
$$U = \bigcup_{\xi \in \Xi(M)} U_\xi $$
be an open neighborhood of $\crit(f) \subset M$, so that 
\begin{enumerate}
\item $\alpha$ is non-zero on $M - U$;
\item each $U_\xi$ is an open neighborhood of $C_\xi$; 
\item $U_\xi$ are pairwise disjoint; and
\item each $U_\xi$ intersects only orbit-type
strata whose closures intersect $C_\xi$.
\end{enumerate}
Any union of sufficiently small neighborhoods $U_\xi$ of $C_\xi$
has these properties.  Since $X_f$ is tangent to the $K$-orbits,
$\phi$ is non-zero on $M - U$.  By \eqref{pos}, for any $R > 0$
there exists an $s(R)$ such that
$ \Vert \Phi(m) + s \phi(m) \Vert > R $
for all $m \in M - U $ and $s > s(R)$.  Since $h$ has support in a
ball of some radius $R(h)$,
\begin{equation} \label{rapid}
  \int_{M - U} \sum_I \eta_I \wedge \exp(\omega + s \d \alpha) 
\partial_I h(\Phi(m) + s \phi(m)) = 0
\end{equation}
for $s > s(R(h))$.  Let $\mu_{\xi,s} := \mu_{U_\xi,K,s}(\eta)$ denote
the distribution defined \eqref{mus} except that integration is over
$U_\xi$.  By \eqref{rapid}, for $s$ sufficiently large
\begin{equation}\label{musum} (\mu_{M,K}(\eta),h) = 
\sum_{\xi \in \Xi(M)} (\mu_{\xi,s},h) .
\end{equation}

\subsection{The limit distributions}

We will show that $\mu_{\xi,s}$ has a distributional limit
$\mu_{\xi,\infty}$ as $s \to \infty$.  In most of this section we will
assume that $\xi$ is central, that is, fixed by the coadjoint action
of $K$.

\begin{lemma}  \label{raycor} \label{genred} (compare Paradan
\cite[3.8]{pa:mo})
\begin{enumerate}
\item Let $[s_1,s_2]$ and $\lambda \in \k^*$ be such that for all $s
\in [s_1,s_2]$, $\Phi_s(\partial U_\xi)$ does not contain $\lambda$.
The restriction of $\mu_{\xi,s}$ to a neighborhood of $\lambda$ is
independent of $s \in [s_1,s_2]$.
\item 
$\mu_{\xi,s}$ converges to a limit $\mu_{\xi,\infty}$ as $s \to
\infty$.
\end{enumerate}
\end{lemma}

\begin{proof} 
(a) follows from \ref{bound}.  (b) $(\phi,\phi)$ is bounded from below
by a positive constant on $\partial U_\xi$.  By \eqref{pos2},
$(\Phi_s,\Phi_s)$ is bounded from below by a constant that approaches
infinity as $s$ does.  Hence for any $\lambda \in \k^*$, there exists
an $s > 0$ such that $\Phi_{s_1}(\partial U_\xi)$ does not contain
$\lambda$ for $s_1 > s$.  The claim follows from (a).
\end{proof} 

We say that $\xi \in \Xi(M)$ is {\em minimal} if $(\xi,\xi)$ is the
minimum value of $(\Phi,\Phi)$.

\begin{lemma}  \label{regval}  Suppose that $\xi$ is minimal.
\begin{enumerate}
\item $\Phi_s(\partial U_\xi)$ does not contain $\xi$ for any $s \in
[0,\infty)$;
\item $\mu_{\xi,\infty}$ is equal to $\mu_{\xi,0}$ in a neighborhood
of $\xi$;
\item $\mu_{\xi,\infty}$ is polynomial on any ray beginning at $\xi$.
\end{enumerate}
\end{lemma}

\begin{proof}  
(a) $(\xi,\xi)$ is the minimum of $(\Phi,\Phi)$ and $C_\xi =
\Phinv(K\xi)$, so $(\Phi,\Phi) > (\xi,\xi)$ on $\partial U_\xi$.
Hence $(\Phi_s,\Phi_s) = (\Phi,\Phi) + 2s(\Phi,\phi) + s^2 (\phi,\phi)
\ge (\Phi,\Phi) > (\xi,\xi)$ on $\partial U_\xi$.  (b) follows from
(c) and \ref{genred} (a).  (c) The fixed point set $U_\xi^{\lambda -
\xi}$ of $\lambda - \xi$ satisfies $ (\Phi_s(U_\xi^{\lambda -
\xi}),\lambda - \xi) = (\xi,\lambda - \xi)$.  Since $(\lambda,\lambda
- \xi) \neq (\xi,\lambda - \xi)$, $\lambda - \xi$ acts infinitesimally
freely on $(\Phi_s,\lambda - \xi)^{-1} ((\lambda,\lambda - \xi))$.
The result now follows from Proposition \ref{ray1}.
\end{proof}

Suppose that $\xi$ is not minimal, and let $U_\xi'$ denote the
Hamiltonian $K_\xi$-manifold obtained from flipping the negative
weights, as in \eqref{vb}, with $ \omega'$ and $\Phi' | U_\xi'$ the
new two-form and moment map.  Since $(\Phi',\xi) \ge (\xi,\xi)$, $\xi$
is minimal for $U_\xi'$.  

\begin{lemma} $(\phi,\xi)$ is non-negative in a neighborhood of $C_\xi$.
\end{lemma}

\begin{proof}   We have $ (\phi(x),\xi) = g(X_f(x),\xi_M(x))$ which
equals $\Vert \xi_M(x) \Vert^2 $ plus a function whose second
derivative vanishes at $C_\xi$.  (One can write the function
explicitly using the local model \eqref{model}).  The function $\Vert
\xi_M \Vert^2$ is Morse-Bott along $Z_\xi$, since $Z_\xi$ is a
component of the fixed point set of $\xi_M$ and $K$ is compact.  The
Hessian of $(\phi,\xi)$ along $Z_\xi$ at $C_\xi$ is equal to the
Hessian of $\Vert \xi_M \Vert^2$, which is positive.  Hence, the
Hessian of $(\phi,\xi)$ along $Z_\xi$ is non-degenerate in a
neighborhood of $C_\xi$ in $Z_\xi$.  Also, $(\phi,\xi)$ is constant on
$Z_\xi$.  It follows that $(\phi,\xi)$ is Morse-Bott along $Z_\xi$, in
a neighborhood of $Z_\xi \cap \Phinv(\xi)$.  Since $(\phi,\xi)$
vanishes on $Z_\xi$ and the Hessian is non-negative, $(\phi,\xi)$ is
non-negative in a neighborhood of $C_\xi$.
\end{proof} 

After shrinking $U_\xi$ if necessary, we may assume that $(\phi,\xi)$
 and $(\phi',\xi)$ are non-negative on $U_\xi$.  Let $\Phi'_s = \Phi'
 + s \phi'$ denote the Witten deformation for $U_\xi'$.  Define
$$ \Phi^u_s | U_\xi = (1-u) \Phi_s | U_\xi + u \Phi'_s | U_\xi .$$
and $ \mu_{\xi,s}^u(\eta)$ the twisted Duistermaat-Heckman
distribution for $\omega_s^u $, as above.

\begin{proposition}  \label{indep} For any $\lambda \in \k^*$, 
for $s$ sufficiently large, $\mu_{\xi,s}^u$ is independent of $u \in
[0,1]$ in a neighborhood of $\lambda$.
\end{proposition} 

\begin{proof}  By \ref{bound} it suffices to show that 
$\Phi^u_s(\partial U_\xi)$ does not contain $\lambda$ for all $u$ and
$s$ sufficiently large.  Since $\phi = \phi'$ on $Z_\xi$ and is
non-vanishing on $Z_\xi \cap \partial U_\xi$, there exists a
neighborhood $V_\xi$ of $Z_\xi \cap \partial U_\xi$ such that $u \phi
+ (1-u) \phi' \neq 0$ on $V_\xi$ for all $u \in [0,1]$.  Hence for $s$
sufficiently large, $\Phi_s^u(V_\xi)$ does not contain $\lambda$.  On
the other hand, the complement $V_\xi^c$ of $V_\xi$ is compact and so
$(\phi,\xi),(\phi',\xi)$ are bounded below on $V_\xi^c$ by a positive
constant.  It follows that $(\Phi_s^u,\Phi^u)$ is bounded below on
$V_\xi^c$ by a constant which approaches infinity as $s$ does.  Hence
$\Phi_s^u(V_\xi^c)$ cannot contain $\lambda$ either. \end{proof}

In the case $\xi$ is central, \ref{indep} and \ref{regval} (b) imply
$ \mu_{\xi,\infty} = \mu_{U_\xi',K_\xi}(\eta) $
in a neighborhood of $\xi$.  By \eqref{vb}
$ \mu_{U_\xi',K_\xi}(\eta) = \mu_{Z_\xi,K_\xi}(\eta \wedge
\Eul(\nu_\xi)^{-1}_\xi) $
in a neighborhood of $\xi$.  This completes the proof of Theorem
\ref{Hform}.

\begin{corollary} \label{tempered}
If $M$ is a Hamiltonian $K$-manifold with proper moment map and finite
number of orbit-type strata, and $\eta \in \Omega_K(M)$ closed, then
$\mu_{M,K}(\eta)$ is a tempered distribution.
\end{corollary}

\begin{proof}  If $M$ has a finite
number of orbit-type strata, then the sum in Theorem \ref{Hform} is
finite.  $\mu_{\xi,\infty}$ is tempered for all $\xi \in \Xi(M)$,
hence $\mu_{M,K}(\eta)$ is a finite sum of tempered distributions.
\end{proof}

\subsection{Further comments}

\begin{enumerate}
\item 
One-parameter localization \ref{Hloc} for central, generic
one-parameter subgroups is a special case of localization via the
norm-square \ref{Hform}.  Indeed, Let $\z \subset \k$ denote the
center of $\k$, and suppose that $\zeta \in \z$. We can use $\zeta$ to
shift the moment map $ \Phi_s = \Phi + s \zeta .$ For sufficiently
large $s$, an element $m \in M$ is fixed by $\Phi_s(m)$ if and only if
it is fixed by $\zeta$.  The subsets $Z_\xi$ are components of
$M^\zeta$, and \ref{Hform} reduces to \ref{Hloc}.

\item The statement and proof of \eqref{Hform} are 
the same in the case that $M$ is a Hamiltonian $K$-orbifold with
proper moment map.
\end{enumerate}

\section{Pairing with invariant functions}
\label{pairsec}

Let $M$ be a Hamiltonian $K$-manifold with proper moment map, and
$\eta \in \Omega_K(M)$ closed.  By \eqref{Ind}, for each $\xi \in
\Xi(M)$, the contribution from $\xi$ to $(\mu_{M,K}(\eta),h)$ is
\begin{equation} \label{contrib}
 ( \mu_{Z^\circ_\xi,K_\xi} ( \Res^{M,K}_{Z^\circ_\xi,K_\xi} \eta
\wedge \Eul(\nu_\xi)^{-1}_\xi ) ,
\Vol_{K_\xi}^K \Res^{K}_{K_\xi} h) .\end{equation}
Hence 
\begin{equation} \label{Hpair}  (\mu_{M,K}(\eta),h) = \sum_{\xi \in \Xi(M)}
(\mu_{Z_\xi^\circ,K_\xi}(\eta \wedge \Eul^{-1}_\xi(\nu_\xi)),
\Vol_{K_\xi}^K \Res^K_{K_\xi} h) .\end{equation}
Suppose that $\Phinv(\xi)$ is contained in the principal orbit-type
stratum for the action of $K_\xi$ on $Z_\xi$.  In this section we show
that the contribution from $\xi$ can be expressed as an integral over
the symplectic quotient
$$Z_{(\xi)} := (Z_\xi)_{(\xi)} = C_\xi/K .$$  
Let $K_\xi'$ denote the identity component of the generic stabilizer
of $K_\xi$ on $\Phinv(\xi) \cap Z_\xi$.  The assumption that
$\Phinv(\xi) \cap Z_\xi$ is contained in the principal orbit-type
stratum of $Z_\xi$ implies that $K'_\xi$ acts trivially on the
annihilator of $\k'_\xi$.  It follows that $K'_\xi$ is normal and the
quotient $K''_\xi := K_\xi/K'_\xi$ is a compact connected Lie group.
Let $K''_{\xi,Z_\xi}$ denote the (finite) generic stabilizer of
$K''_\xi$ on $Z_\xi$.  Let $ \kappa_\xi $ denote the composition of
the restriction $H_K(M) \to H_{K_\xi}(\Phinv(\xi))$ with the
isomorphism
$$ H_{K_\xi}(\Phinv(\xi) \cap Z_\xi) \to H_{K'_\xi}(Z_{(\xi)}) = H(Z_{(\xi)})
\otimes S(\k'^{,*}_\xi)^{K'_\xi} .$$
This extends to forms with smooth coefficients.  In particular, any $h
\in \S(\k_\xi^*)^{K_\xi}$ defines a characteristic class $\kappa_\xi(h)
\in H_{K'_\xi}(Z_{(\xi)})$.  We denote by
$$\nu_{(\xi)} := K'_\xi \backslash ( \nu_\xi |_{\Phinv(\xi)}) \to
Z_{(\xi)}$$ 
the quotient bundle.   

\begin{theorem} \label{Hsmooth}  
If $\Phinv(\xi) \cap Z_\xi$ is contained in the principal orbit-type
stratum in $Z_\xi$, then \eqref{contrib} is equal to
$$ \Vol(K''_\xi/K''_{\xi,Z_\xi}) \int_{Z_{(\xi)} \times \k_\xi'}
\L_{(\xi)} \wedge \Eul(\nu_{(\xi)})_\xi^{-1} \wedge \kappa_\xi(\eta
\wedge \Vol_{K_\xi}^K \Res^K_{K_\xi} h).
 $$
\end{theorem}
Note that the inverted Euler class has tempered-distributional
coefficients, while the remaining factor has coefficients in the ring
of Schwartz functions.  The integral over $\k_\xi'$ refers to the
pairing of these coefficient rings.

\begin{proof} 
Suppose that $\xi = 0 $, $\eta = 1$, $M$ is maximal rank, and
$\Phinv(0)$ is contained in the principal orbit-type stratum for $M$.
Let $\alpha \in \Omega^1(\Phinv(0),\k)^K$ denote a connection one-form
for the action of $K$ on $\Phinv(0)$.  Let $\pi_0 : \ \Phinv(0) \to
M_{(0)}$ denote the projection.  By the coisotropic embedding theorem,
a neighborhood of $\Phinv(0)$ is $K$-symplectomorphic to a
neighborhood of $\Phinv(0) \times \{ 0 \}$ in the Hamiltonian
$K$-manifold
\begin{equation} \label{coiso}
( \Phinv(0) \times \k^*, \ \ \pi_0^*\omega_{(0)} + \d (
\lambda,\alpha ) ) .\end{equation}
%
%
Let $\pi_1,\pi_2$ denote the projections 
$$\Phinv(0)
\stackrel{\pi_1}{\to} T \backslash \Phinv(0) \stackrel{\pi_2}{\to}
M_{(0)} .$$ 
The two-form $\d \alpha$ is $T$-basic and descends to a closed
$\k$-valued two-form $\pi_{1,*} \d \alpha $.  By \eqref{coiso} and
\eqref{Hred}, the volumes of quotients at generic $\lambda \in \k^*$
are
$$ p(\lambda) := \Vol(M_{(\lambda)}) = I_{T \backslash \Phinv(0)} (
\exp(\pi_2^*\omega_{(0)} + ( \lambda , \pi_{1,*} \d \alpha) )).$$
Let $h \in \S(\k^*)^K$, $f = \F_\t^{-1}(\Pi \Res^{\k^*}_{\t^*} h)$ and
  $p = \sum_I p_I \lambda_I$.  Using \eqref{Rvol}
\begin{eqnarray*} 
\Vol(K/K_M)^{-1} (\mu_{M,K},h) &=&  ( (\Vol_T^K)^{-1} p,  h )_{\k^*}  \\
&=&  (\# W)^{-1}  ( \Res_{\t^*}^{\k^*} p,  \Pi h )_{\t^*}  \\
               &=&  (\# W)^{-1} ( \F_\t^{-1} (\Res_{\t^*}^{\k^*} p),
               \F_\t^{-1}(\Pi \Res_{\t^*}^{\k^*}h))_\t \\
               &=& (\# W)^{-1} \sum_I p_I (\partial_I f)(0) \\ 
               &=& (\# W)^{-1} I_{T \backslash \Phinv(0)} (
\exp(\pi_2^*\omega_{(0)}) f(\pi_{1,*} \frac{\d \alpha}{2 \pi i})) \end{eqnarray*}
defined using the formal power series of $f$ at $0$.  Choose an
orthonormal basis $\xi_1,\ldots,\xi_r$ for $\t$ and
$\xi_{r+1},\ldots,\xi_n$ for $\k/\t \cong \t^\perp$.  We can replace
the integral over $T \backslash \Phinv(0)$ with
$$ \frac{1}{\Vol(T)} \int_{\Phinv(0)} ( \exp(\pi_2^* \omega_{(0)})
\wedge f(\frac{\d \alpha}{2 \pi i}) \wedge \prod_{j=1}^r (\alpha,\xi_j)) .$$
It remains to integrate over the fiber $K$ of $\Phinv(0) \to M_{(0)}$.
Writing 
$$ d\alpha = \pi_0^* \curv(\alpha) - \hh [\alpha, \alpha] $$
we see that the component of 
$ f(\frac{\d \alpha}{2 \pi i}) \wedge \prod_{j=1}^r (\alpha,\xi_j)
$
that contributes to the fiber integral is 
$$ ( \F_t^{-1}\Pi f) \left( \pi_0^*
\frac{\curv(\alpha)}{2\pi i} ) \right) \wedge 
\prod_{i=1}^n
(\alpha,\xi_j) .$$
%
%
%
Integrating over the fiber changes the factor $\prod_{i=1}^n
(\alpha,\xi_j) $ to $\Vol(K)$.  We have
\begin{eqnarray*}
 (\# W)^{-1} (\F_t^{-1}(\Pi^2 \Res^{\k^*}_{\t^*}
h))( \frac{\curv(\alpha)}{2\pi i}) &=& (\F_\k^{-1}(h))(
\frac{\curv(\alpha)}{2\pi i}) \\
& =&
    \L_0 \wedge \kappa_0(h) .\end{eqnarray*}
It follows that 
$$ (\mu_{M,K},h) =  \Vol(K/K_M) I_{M_{(0)}} (\L_0
\wedge \kappa_0(h)) .$$
The general case is the same, except that the forms on $Z_{(\xi)}$ are
$K_\xi'$-equivariant and the form $\eta$ is to be included.
\end{proof}

\begin{corollary}  Let $M$ be a compact Hamiltonian $K$-manifold,
$\eta \in \Omega_K(M)$ closed, and $h \in \S(\k^*)^K$.  Suppose that
$\Phinv(\xi) \cap Z_\xi$ is contained in the principal orbit-type
stratum of $Z_\xi$, for all $\xi \in \Xi(M)$.  The pairing
$(I_{M,K}(\eta),\F_\k^{-1}(h))$ is equal to
$$ \sum_{\xi \in \Xi(M)} 
\Vol(K''_\xi/K''_{\xi,Z_\xi}) \int_{Z_{(\xi)} \times \k_\xi'}
 \kappa_\xi(\eta
\wedge \Vol_{K_\xi}^K \Res^K_{K_\xi} h) \wedge
\Eul(\nu_{(\xi)})_\xi^{-1}.
$$
\end{corollary}

\begin{proof}  Consider the family of equivariant symplectic forms
$\eps \omega_K$ for $\eps \in (0,1]$.  The corollary follows from
taking the limit $\eps \to 0$ of \eqref{Hpair}, using Theorem
\ref{Hsmooth}: The stratification is independent of $\eps$, and
$\L_{(\xi)} \to 1$ as $\eps \to 0$.
\end{proof}

If $\Phinv(\xi) \cap Z_\xi$ is not contained in the principal
orbit-type stratum of $Z_\xi$, then \eqref{contrib} can be written as
a finite sum of integrals over symplectic quotients near $\xi$, by the
gluing rule in Meinrenken \cite{me:sym}, but I do not know a nice
formula for the limit.


\section{Examples} 
\label{examples}

In these examples we will compare the one-parameter and norm-square
localization formulas.

\subsection{ $K = U(1)$ acting on $M = \P^1$} 
\begin{figure}[h]
\setlength{\unitlength}{0.00053333in}
\begingroup\makeatletter\ifx\SetFigFont\undefined%
\gdef\SetFigFont#1#2#3#4#5{%
  \reset@font\fontsize{#1}{#2pt}%
  \fontfamily{#3}\fontseries{#4}\fontshape{#5}%
  \selectfont}%
\fi\endgroup%
{\renewcommand{\dashlinestretch}{30}
\begin{picture}(4222,1364)(0,-10)
\thicklines
\path(3000,352)(600,352)
\path(720.000,382.000)(600.000,352.000)(720.000,322.000)
\path(720.000,82.000)(600.000,52.000)(720.000,22.000)
\path(600,52)(1800,52)
\put(-300,252){\makebox(0,0)[lb]{{\SetFigFont{12}{14.4}{\rmdefault}{\mddefault}{\updefault}equals}}}
\put(-300,-48){\makebox(0,0)[lb]{{\SetFigFont{12}{14.4}{\rmdefault}{\mddefault}{\updefault}minus}}}
\path(1800,1027)(4200,1027)
\path(4080.000,997.000)(4200.000,1027.000)(4080.000,1057.000)
\path(3000,727)(4200,727)
\path(4080.000,697.000)(4200.000,727.000)(4080.000,757.000)
\put(900,927){\makebox(0,0)[lb]{{\SetFigFont{12}{14.4}{\rmdefault}{\mddefault}{\updefault}equals}}}
\put(900,627){\makebox(0,0)[lb]{{\SetFigFont{12}{14.4}{\rmdefault}{\mddefault}{\updefault}minus}}}
\path(1800,1327)(3075,1327)
\end{picture}
}
\caption{One-parameter localization for $\P^1$ \label{p1ex1fig}}
\end{figure}
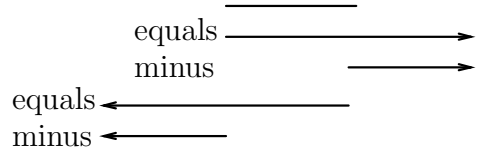

 We identify the Lie algebra $i\R$ of $U(1)$ with $\R$ by division by
$2\pi i$.  If we choose on $\R$ the standard inner product then the
weight lattice becomes identified with $\Z$, $\mu_{\k^*}$ is Lebesgue
measure on $\R$, and the volume of $K$ is $1$.  The action of $U(1)$
on $\P^1$ by $ z[w_0,w_1] = [z^{-a} w_0,z^{-b} w_1] $ has moment map
$$ \Phi([w_0,w_1]) = \frac{ a |w_0|^2 + b |w_1|^2}{|w_0|^2 + |w_1|^2} .$$
There are two fixed points, at $w_0 = 0$, resp. $w_1 = 0$.  The
tangent weights at the fixed points are $\pm (b-a)$.  Let $H_\pm$
denote the Heaviside distributions, equal to $\mu_{\k^*}$ on the
positive (resp. negative) real numbers and zero elsewhere.  
One-parameter localization with $\xi > 0$ gives
$$ \mu_{M,K} = \frac{1} {b-a} (\delta(a)\star H_+ - \delta(b) \star
H_+) = \frac{\chi_{[a,b]} \mu_{\k^*} }{b-a} $$
where $\chi_{[a,b]}$ the characteristic function for the interval
$[a,b]$ and $\star$ denotes convolution.  For negative action chamber
$\xi < 0$ localization gives
$$ \mu_{M,K} = \frac{1} {b-a} (- \delta(a) \star H_- + \delta(b) \star
H_-) = \frac{\chi_{[a,b]}  \mu_{\k^*}}{b-a} .$$
This is shown graphically in Figure \ref{p1ex1fig}.

%

The Kirwan-Ness localization is as follows for $a < 0 < b $ is as
follows.  The critical components of $f$ are the two $T$-fixed points
and the zero level set:
$$  C_{a} = \{ w_0 = 0 \} = M_{a}, \ \ C_{b} = \{ w_1 = 0 \} = M_{b} $$
$$ C_0 = \Phinv(0), \ M_0 = \{ [w_1,w_2], w_1w_2 \neq 0 \} = M - M_{a}
- M_{b}.$$

\begin{figure}[h]
\setlength{\unitlength}{0.00033333in}
\begingroup\makeatletter\ifx\SetFigFont\undefined%
\gdef\SetFigFont#1#2#3#4#5{%
  \reset@font\fontsize{#1}{#2pt}%
  \fontfamily{#3}\fontseries{#4}\fontshape{#5}%
  \selectfont}%
\fi\endgroup%
{\renewcommand{\dashlinestretch}{30}
\begin{picture}(5797,1289)(0,-10)
\thicklines
\texture{44555555 55aaaaaa aa555555 55aaaaaa aa555555 55aaaaaa aa555555 55aaaaaa 
	aa555555 55aaaaaa aa555555 55aaaaaa aa555555 55aaaaaa aa555555 55aaaaaa 
	aa555555 55aaaaaa aa555555 55aaaaaa aa555555 55aaaaaa aa555555 55aaaaaa 
	aa555555 55aaaaaa aa555555 55aaaaaa aa555555 55aaaaaa aa555555 55aaaaaa }
\path(3075,1252)(4275,1252)
\path(3075,1252)(4275,1252)
\path(1695.000,682.000)(1575.000,652.000)(1695.000,622.000)
\path(1575,652)(5775,652)
\path(1575,652)(5775,652)
\path(5655.000,622.000)(5775.000,652.000)(5655.000,682.000)
\path(1695.000,82.000)(1575.000,52.000)(1695.000,22.000)
\path(1575,52)(3075,52)
\path(1575,52)(3075,52)
\path(4275,52)(5775,52)
\path(4275,52)(5775,52)
\path(5655.000,22.000)(5775.000,52.000)(5655.000,82.000)
\put(3675,52){\makebox(0,0)[lb]{{\SetFigFont{12}{14.4}{\rmdefault}{\mddefault}{\updefault}}}}
\put(0,577){\makebox(0,0)[lb]{{\SetFigFont{12}{14.4}{\rmdefault}{\mddefault}{\updefault}equals}}}
\put(0,52){\makebox(0,0)[lb]{{\SetFigFont{12}{14.4}{\rmdefault}{\mddefault}{\updefault}minus}}}
\end{picture}
}
\caption{Norm-square localization for $\P^1$ \label{p1ex5fig}}
\end{figure}
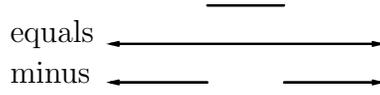

For $\xi = 0 $, $Z_\xi^\circ = M_\xi$ is the complement of the
$K$-fixed points.  The Duistermaat-Heckman measure for $Z_\xi^\circ$
is $ \mu_{\k^*} \chi_{[a,b]}/(b-a)$.  Its unique extension which is
piecewise polynomial on any ray beginning at $0$ is
$\mu_{\k^*}/(b-a)$.  Theorem \ref{Hform} gives
\begin{eqnarray*} \mu_{M,K} &=& \mu_{M,0} + \mu_{M,b} + \mu_{M,a} \\
&=& \frac{1}{b-a} \left( \mu_{\k^*}- H_- \star \delta(a) - 
H_+ \star \delta(b) \right)  \\
&=& \frac{1}{b-a}  \mu_{\k^*}(1 - \chi_{(-\infty,a]} - \chi_{[b,\infty)}).
\end{eqnarray*}
The formula is shown graphically in Figure \ref{p1ex5fig}.

\subsection{ $SU(3)$ acting on a $G_2$-coadjoint orbit}

\begin{figure}[h]
\setlength{\unitlength}{0.00013333in}
\begingroup\makeatletter\ifx\SetFigFont\undefined%
\gdef\SetFigFont#1#2#3#4#5{%
  \reset@font\fontsize{#1}{#2pt}%
  \fontfamily{#3}\fontseries{#4}\fontshape{#5}%
  \selectfont}%
\fi\endgroup%
{\renewcommand{\dashlinestretch}{30}
\begin{picture}(8811,9055)(0,-10)
\path(3343,8112)(3343,6312)
\blacken\path(3328.000,6372.000)(3343.000,6312.000)(3358.000,6372.000)(3343.000,6390.000)(3328.000,6372.000)
\path(6971,8128)(5412,7228)
\blacken\path(5456.463,7270.988)(5412.000,7228.000)(5471.462,7245.007)(5479.552,7266.997)(5456.463,7270.988)
\path(6971,8128)(6971,6328)
\blacken\path(6956.000,6388.000)(6971.000,6328.000)(6986.000,6388.000)(6971.000,6406.000)(6956.000,6388.000)
\path(8799,4994)(7240,5894)
\blacken\path(7299.462,5876.993)(7240.000,5894.000)(7284.463,5851.012)(7307.552,5855.003)(7299.462,5876.993)
\path(8799,4994)(7240,4094)
\blacken\path(7284.463,4136.988)(7240.000,4094.000)(7299.462,4111.007)(7307.552,4132.997)(7284.463,4136.988)
\path(6971,1812)(6971,12)
\blacken\path(6956.000,72.000)(6971.000,12.000)(6986.000,72.000)(6971.000,90.000)(6956.000,72.000)
\path(3296,1828)(3296,29)
\blacken\path(3281.000,89.000)(3296.000,29.000)(3311.000,89.000)(3296.000,107.000)(3281.000,89.000)
\path(3318,1796)(1759,896)
\blacken\path(1803.463,938.988)(1759.000,896.000)(1818.462,913.007)(1826.552,934.997)(1803.463,938.988)
\path(3296,8128)(1737,9028)
\blacken\path(1796.462,9010.993)(1737.000,9028.000)(1781.463,8985.012)(1804.552,8989.003)(1796.462,9010.993)
\path(1618,4962)(59,4062)
\blacken\path(103.463,4104.988)(59.000,4062.000)(118.462,4079.007)(126.552,4100.997)(103.463,4104.988)
\path(1571,4978)(12,5878)
\blacken\path(71.462,5860.993)(12.000,5878.000)(56.463,5835.012)(79.552,5839.003)(71.462,5860.993)
\path(6980,1811)(5421,2711)
\blacken\path(5480.462,2693.993)(5421.000,2711.000)(5465.463,2668.012)(5488.552,2672.003)(5480.462,2693.993)
\put(6365,7291){\makebox(0,0)[lb]{{\SetFigFont{20}{24.0}{\rmdefault}{\mddefault}{\updefault}+}}}
\put(2837,7719){\makebox(0,0)[lb]{{\SetFigFont{20}{24.0}{\rmdefault}{\mddefault}{\updefault}-}}}
\put(7749,4934){\makebox(0,0)[lb]{{\SetFigFont{20}{24.0}{\rmdefault}{\mddefault}{\updefault}-}}}
\put(6282,1225){\makebox(0,0)[lb]{{\SetFigFont{20}{24.0}{\rmdefault}{\mddefault}{\updefault}+}}}
\put(2738,813){\makebox(0,0)[lb]{{\SetFigFont{20}{24.0}{\rmdefault}{\mddefault}{\updefault}-}}}
\put(612,4950){\makebox(0,0)[lb]{{\SetFigFont{20}{24.0}{\rmdefault}{\mddefault}{\updefault}+}}}
\end{picture}
}
\hskip 1in
\setlength{\unitlength}{0.00023333in}
\begingroup\makeatletter\ifx\SetFigFont\undefined%
\gdef\SetFigFont#1#2#3#4#5{%
  \reset@font\fontsize{#1}{#2pt}%
  \fontfamily{#3}\fontseries{#4}\fontshape{#5}%
  \selectfont}%
\fi\endgroup%
{\renewcommand{\dashlinestretch}{30}
\begin{picture}(5238,4559)(0,-10)
\path(3928,4530)(3928,30)
\path(3925,30)(28,2280)
\path(3928,4530)(31,2280)
\path(5203,2280)(1306,30)
\path(1304,33)(1304,4532)
\path(5203,2280)(1306,4530)
\texture{44000000 aaaaaa aa000000 8a888a 88000000 aaaaaa aa000000 888888 
	88000000 aaaaaa aa000000 8a8a8a 8a000000 aaaaaa aa000000 888888 
	88000000 aaaaaa aa000000 8a888a 88000000 aaaaaa aa000000 888888 
	88000000 aaaaaa aa000000 8a8a8a 8a000000 aaaaaa aa000000 888888 }
\shade\path(3922,4526)(2622,3776)(3932,3016)(3932,4526)
\path(3922,4526)(2622,3776)(3932,3016)(3932,4526)
\shade\path(1312,3026)(12,2276)(1322,1516)(1322,3026)
\path(1312,3026)(12,2276)(1322,1516)(1322,3026)
\shade\path(3922,1536)(2622,786)(3932,26)(3932,1536)
\path(3922,1536)(2622,786)(3932,26)(3932,1536)
\texture{ffffffff ffeeeeee eeffffff fffbfbfb fbffffff ffeeeeee eeffffff ffbfbbbf 
	bbffffff ffeeeeee eeffffff fffbfbfb fbffffff ffeeeeee eeffffff ffbfbfbf 
	bfffffff ffeeeeee eeffffff fffbfbfb fbffffff ffeeeeee eeffffff ffbfbbbf 
	bbffffff ffeeeeee eeffffff fffbfbfb fbffffff ffeeeeee eeffffff ffbfbfbf }
\shade\path(1302,4526)(2612,3766)(1302,3006)(1302,4516)
\path(1302,4526)(2612,3766)(1302,3006)(1302,4516)
\shade\path(3916,3041)(5226,2281)(3916,1521)(3916,3031)
\path(3916,3041)(5226,2281)(3916,1521)(3916,3031)
\shade\path(1314,1532)(2624,772)(1314,12)(1314,1522)
\path(1314,1532)(2624,772)(1314,12)(1314,1522)
\end{picture}
}
\caption{One-parameter localization for $G_2/P$
\label{g2ex1fig}}
\end{figure}
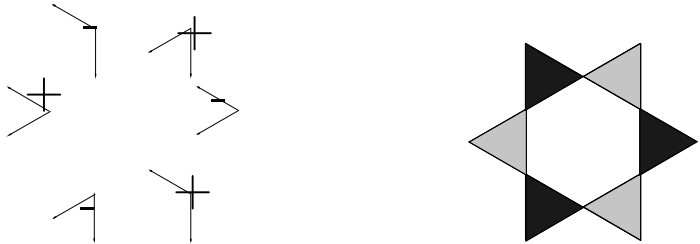

In this example, we apply the localization formulas to the action of
$SU(3) \subset G_2$ on a coadjoint orbit of $G_2$.  Let $K = SU(3)$,
and $\omega_1,\omega_2$ the fundamental weights.  Let $G_2$ denote the
connected simple complex group of type $G_2$.  The dual positive Weyl
chamber for $G_2$ is the span of $\omega_1$ and $\omega_1 + \omega_2$.
Let $P_{\omega_1 + \omega_2}$ denote the maximal parabolic of $G_2$,
so that $M = G_2/P_{\omega_1 + \omega_2}$ is diffeomorphic to the
coadjoint orbit through $\omega_1 + \omega_2$.  The Weyl group $W$ for
$SU(3)$ acts simply transitively on the $T$-fixed points.  First we
compute the Duistermaat-Heckman measure using ordinary localization.
The contribution to $\mu_{M,T}$ from the fixed point $x(w)$
corresponding to $w \in W$ is
$$ \delta_{w \mu} \star \prod_{j =1}^5  \pm H_{ \pm w \beta_j} $$
where $\beta_j$ are the positive roots of $G_2$ not vanishing at
$\omega_1 + \omega_2$
$$ 2 \omega_1 - \omega_2, \  -\omega_1 + 2 \omega_2,\  \omega_1 + \omega_2,
\ 3 \omega_1,\  3 \omega_2 $$  
and the signs are determined by the action chamber.  The $\beta_j$
that are not roots of $SU(3)$ are $ \beta_5 = 3 \omega_1,\ \beta_6 = 3
\omega_2 .$ By \eqref{TtoG}
$$   \mu_{M,K} = 
\frac{1}{6}  \Ind_T^K 
\sum_{w \in W} (-1)^{l(w)} \delta_{w (\omega_1 + \omega_2)}
\star (\pm H_{\pm 3 w\omega_1}) \star (\pm H_{\pm 3 w\omega_2}) .$$
The contributions are shown in Figure \ref{g2ex1fig}.  Each contribution is $\pm
\mu_{\t^*} / (9 \Vert \omega_1 \Vert \Vert \omega_2 \Vert)$ where
non-negative.  Positive (resp. negative) contributions are shown in
light (resp. dark) shading.  The moment polytope for $M$ is
$$ P = \hull(\omega_1,\omega_2,\omega_1 + \omega_2) .$$
Let $F_1$ be the open face connecting $\omega_2,\omega_1 + \omega_2$,
$F_2$ the open face connecting $\omega_1,\omega_1 + \omega_2$, and
$F_3$ the open face connecting $\omega_1,\omega_2$.  Let $F_{ij} = F_1
\cap F_2$.

We compute the Kirwan-Ness stratification.  The inverse image
$\Phinv(F_{12})$ contains a unique point, $x(1) \in M$, which is
$T$-fixed.  None of the other $T$-fixed points map to $\t_+^*$.
Therefore, the remaining points in $\Phinv(\on{int}(\t_+^*))$ have
one-dimensional stabilizers.  Since $\Phinv(\on{int}(\t_+^*))$ has dimension $ 2
\dim(T)$, it is a toric manifold, so the inverse image of any face $F
\subset \on{int}\t_+^*$ has infinitesimal stabilizer the annihilator
of the tangent space of $F$.  The stabilizers of the faces
$F_1,F_2,F_3$ are
$$ \t_1 = \on{span}(h_1), \ \t_2 = \on{span}(h_2), \ \t_3 = \on{span}(h_3) $$
where $h_1,h_2,h_3$ are the coroots of $SU(3)$.  The level set
$\Phinv((\omega_1 + \omega_2)/2)$ is critical with $\xi =
(\omega_1 + \omega_2)/2$.  
\begin{figure}[h]
\setlength{\unitlength}{0.00027489in}
\begingroup\makeatletter\ifx\SetFigFont\undefined
\def\x#1#2#3#4#5#6#7\relax{\def\x{#1#2#3#4#5#6}}%
\expandafter\x\fmtname xxxxxx\relax \def\y{splain}%
\ifx\x\y   
\gdef\SetFigFont#1#2#3{%
  \ifnum #1<17\tiny\else \ifnum #1<20\small\else
  \ifnum #1<24\normalsize\else \ifnum #1<29\large\else
  \ifnum #1<34\Large\else \ifnum #1<41\LARGE\else
     \huge\fi\fi\fi\fi\fi\fi
  \csname #3\endcsname}%
\else
\gdef\SetFigFont#1#2#3{\begingroup
  \count@#1\relax \ifnum 25<\count@\count@25\fi
  \def\x{\endgroup\@setsize\SetFigFont{#2pt}}%
  \expandafter\x
    \csname \romannumeral\the\count@ pt\expandafter\endcsname
    \csname @\romannumeral\the\count@ pt\endcsname
  \csname #3\endcsname}%
\fi
\fi\endgroup
{\renewcommand{\dashlinestretch}{30}
\begin{picture}(16955,5259)(0,-10)
\texture{aaffffff ffaaaaaa aaffffff ffaaaaaa aaffffff ffaaaaaa aaffffff ffaaaaaa 
	aaffffff ffaaaaaa aaffffff ffaaaaaa aaffffff ffaaaaaa aaffffff ffaaaaaa 
	aaffffff ffaaaaaa aaffffff ffaaaaaa aaffffff ffaaaaaa aaffffff ffaaaaaa 
	aaffffff ffaaaaaa aaffffff ffaaaaaa aaffffff ffaaaaaa aaffffff ffaaaaaa }
\shade\path(41,1812)(1616,912)(1616,2712)(41,1812)
\path(41,1812)(1616,912)(1616,2712)(41,1812)
\shade\path(6161,1857)(7736,957)(7736,2757)(6161,1857)
\path(6161,1857)(7736,957)(7736,2757)(6161,1857)
\texture{55888888 88555555 5522a222 a2555555 55888888 88555555 552a2a2a 2a555555 
	55888888 88555555 55a222a2 22555555 55888888 88555555 552a2a2a 2a555555 
	55888888 88555555 5522a222 a2555555 55888888 88555555 552a2a2a 2a555555 
	55888888 88555555 55a222a2 22555555 55888888 88555555 552a2a2a 2a555555 }
\shade\path(7737,4565)(7736,2757)(9318,3680)
	(7736,4580)(7728,2772)
\path(7737,4565)(7736,2757)(9318,3680)
	(7736,4580)(7728,2772)
\path(71.000,5067.000)(41.000,5187.000)(11.000,5067.000)
\path(41,5187)(41,12)
\path(4434.067,2513.029)(4523.000,2599.000)(4404.073,2564.994)
\path(4523,2599)(41,12)
\path(41,1812)(1616,912)
\path(6191.000,5112.000)(6161.000,5232.000)(6131.000,5112.000)
\path(6161,5232)(6161,57)
\path(10554.067,2558.029)(10643.000,2644.000)(10524.073,2609.994)
\path(10643,2644)(6161,57)
\path(7770.834,4450.959)(7741.000,4571.000)(7710.834,4451.042)
\path(7741,4571)(7736,957)
\path(9268.323,3630.765)(9357.000,3717.000)(9238.175,3682.641)
\path(9357,3717)(6160,1859)
\path(6161,1857)(7736,957)
\path(7752,4572)(9322,3676)
\path(12491.000,5112.000)(12461.000,5232.000)(12431.000,5112.000)
\path(12461,5232)(12461,57)
\path(16854.067,2558.029)(16943.000,2644.000)(16824.073,2609.994)
\path(16943,2644)(12461,57)
\path(14047,4576)(15622,3676)
\shade\path(14045,4562)(14045,2762)(15627,3685)
	(14045,4585)(14037,2777)
\path(14045,4562)(14045,2762)(15627,3685)
	(14045,4585)(14037,2777)
\put(5133,2645){\makebox(0,0)[lb]{=}}
\put(11412,2580){\makebox(0,0)[lb]{-}}
\end{picture}
}
\caption{Norm-square localization for $G_2/P$ \label{g2exfig}}
\end{figure}
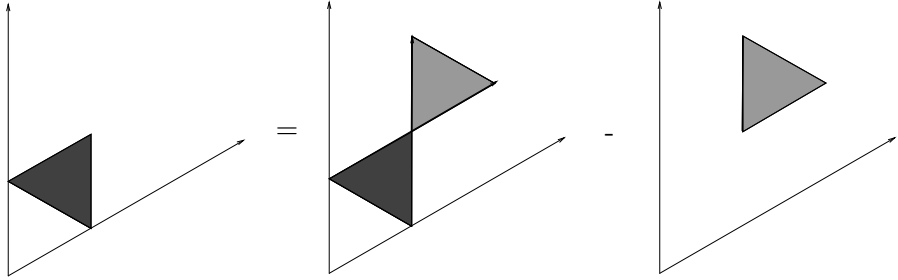
The fixed point component $Z_\xi$ has moment image
$$ \Phi(Z_\xi) = \hull(2\omega_2 - \omega_1, 2\omega_1 - \omega_2 ) .$$
The unstable manifold $Y_\xi$ has image under the moment map for $T$
$$ \on{proj}^\k_\t \Phi(\ol{Y_\xi}) = \hull(2\omega_2 - \omega_1, 2\omega_1 -
\omega_2 , \omega_1 + \omega_2) .$$
None of the other faces $F_j$ contain points $\xi$ with $\xi \in
\t_j$.  Therefore, there are no other critical points in $
\Phinv(\on{int}(\t_+^*))$.  Finally consider the inverse image of the
vertices $F_{jk} = F_{13} $ or $ F_{23}$.  $\Phinv(F_{jk})$ does not
contain a $T$-fixed point.  $\Phinv(F_{jk})$ does not contain a point
$m$ stabilized by $\on{span}(F_{jk})$.  Indeed, since the stabilizer
$K_m$ does not contain a maximal torus, $K_m$ cannot intersect the
semisimple part $[K_{\Phi(m)}, K_{\Phi(m)}]$.  Therefore, $K_m$ is
one-dimensional.  Let $X$ denote the fixed point component of $K_m$
containing $m$.  Since $K_m$ is one-dimensional, the image $\Phi(X)$
is codimension one, and so meets $\Phinv(\on{int}(t_+^*))$.  This
implies that the $\k_m$ is conjugate to either $\t_j$ or $\t_k$, and
so $\k_m$ cannot equal the span of $F_{jk}$.  Therefore,
$$ \Xi(M) = \{ \omega_1 + \omega_2, \hh (\omega_1 + \omega_2) \} .$$
One can show that the Kirwan-Ness stratification coincides with the orbit
stratification for $G$, just as in the previous example.  In
particular $M$ is a two-orbit variety, with one open orbit and one of
complex codimension two.

The contributions to the norm-square localization formula can be
described as follows.  For $\xi = (\omega_1 + \omega_2) /2 $ we have
$$ \mu_{\xi,\infty} =  \delta( (\omega_1 + \omega_2,\xi)
= \hh \Vert \omega_1 + \omega_2 \Vert) \star
H_{\omega_1 + \omega_2} / (9 \Vert \omega_1 \Vert
\Vert \omega_2 \Vert) .$$
where $H_{\omega_1 + \omega_2}$ is the Heaviside distribution for
$(\omega_1 + \omega_2 , \xi) \ge 0$.  Therefore,
$$ \Ind_{T}^K \mu_\xi^\infty = \Ind_T^K (\chi_{P} -
\chi_{Q} ) \mu_{\t^*} / (9 \Vert
\omega_1 \Vert \Vert \omega_2 \Vert) .$$
where $\chi_{P}, \chi_{Q}$ are the characteristic functions
for the polytope $P$, resp. the cone
$$ Q = \R_{\le 0} (P - (\omega_1 + \omega_2)) + \omega_1 + \omega_2 .$$
For $\xi = \omega_1 + \omega_2$ we get
$$ \Ind_{T}^K \mu_\xi^\infty = \Ind_T^K \frac{ \chi_{Q}  \mu_{\t^*}}{9 \Vert \omega_1 \Vert \Vert \omega_2
\Vert} .$$
Hence
$$ \mu_{M,K} = \Ind_T^K \left( \frac{\chi_{P} \mu_{\t^*}}
{9 \Vert \omega_1
\Vert \Vert \omega_2 \Vert}
\right) \hfill .$$
See Figure \ref{g2exfig}.

\section{A remark on sheaf cohomology}
\label{alg}

In algebraic geometry there is a formula which expresses the index of
a sheaf of a stratified variety as a sum over the strata.  Let $G$ be
a reductive complex group, $R(G)$ the ring of finite linear
combinations of irreducible characters, and ${\RR(G)} = \Hom(R(G),\Z)$
its dual.  Let $M$ be a smooth quasiprojective variety, and $E \to M$
a $G$-equivariant vector bundle.  The {\em equivariant index} of $E$
is the virtual representation
$$ I_{M,K}(E) = \sum_{j = 0}^{\dim(M)} (-1)^j H^j(M,E) .$$
We will assume that the multiplicity of any irreducible representation
is finite, so that $I_{M,K}(E)$ defines an element in ${\RR(G)}$. 

Suppose $M$ decomposes into a disjoint union of smooth $G$-stable
subvarieties
$$ M = \bigcup_{\xi \in \Xi(M)} M_\xi .$$
Let $T_{M_\xi} M$ is the normal bundle of $M_\xi \to M$, $T_{M_\xi}^*
M$ its dual.   The Euler class
$$ \Eul(T_{M_\xi} M) := \Lambda^{\on{even}}(T_{M_\xi} M)
\ominus \Lambda^{\on{odd}}(T_{M_\xi} M) $$
has a formal inverse
$$ \Eul(T_{M_\xi}M)^{-1} := (-1)^{\codim(M_\xi)} \det(T_{M_\xi}^* M)
\otimes S(T_{M_\xi}^* M) $$
where $S$ resp. $\Lambda$ denotes the direct sum of symmetric
resp. exterior powers and $\det$ the top exterior power.  Let
$\Res^M_{M_\xi}$ denote restriction to $M_\xi$.  The
Cousin-Grothendieck spectral sequence (take the Euler characteristic
of the local cohomologies) produces a formula in ${\RR(G)}$
\begin{equation} \label{cr}
 I_{M,K}(E) = \sum_{\xi \in \Xi(M)} I_{M_\xi,K}( \Res^M_{M_\xi} E
\otimes \Eul(T_{M_\xi}M )^{-1}) \end{equation}
assuming that the representations on the right hand side have finite
multiplicities, see Teleman \cite{te:qu} and Hartshorne \cite[Section
4]{ha:rd}.

In some sense the localization theorems in equivariant cohomology or
K-theory are attempts to extend this result to manifolds with group
actions; so far this has been done only in special cases.  One
parameter localization arises from the stratification defined by the
action of a circle subgroup $G = \C^*$.  Let $\Xi(M)$ denote the set
of connected components of the fixed point set $M^G$ in $M$.  For any
$\xi \in \Xi(M)$, let
$ M_\xi = \{ m \in M, \ \ \lim_{z \to 0} zm \in \xi \} $
denote the stable manifolds for the flow generated by the action, as
in Bialinicki-Birula \cite{bi:so}.  The formula \eqref{cr} is a
sheaf-theoretic version of \ref{Hloc}.  \eqref{cr} applied to the
Kirwan-Ness stratification gives a sheaf-theoretic version of
\ref{Hform}.  It remains an open question, at least for me, whether
there is a more general formula in equivariant $K$-theory or in
equivariant de Rham theory analogous to \eqref{cr}.  For instance, the
decomposition of a projective spherical $G$-variety into $G$-orbits
produces a formula in K-theory not covered by one-parameter
localization or localization via the norm-square of the moment map.

\section{2d Yang-Mills}

The basic reference for mathematical two dimensional Yang-Mills theory
is Atiyah-Bott \cite{at:mo}.  Let $K$ denote a connected compact Lie
group and $G$ the complexification of $K$.  Fix the basic inner
product $(\ , \ ): \k \times \k \to \R$ and use it to identify $\k$
with its dual $\k^*$.  Let $P$ be a principal $K$-bundle and $P(\k) =
P \times_K \k$ the adjoint bundle.  Similarly, let $P(G) = P \times_K
G$ the associated principal $G$-bundle.  Let $\Omega^\bullet(X,P(\k))$
the space of forms with values in $P(\k)$.  The inner product on $\k$
induces a metric on $P(\k)$.  Combining this with the wedge product
gives map
$$ \Omega^k(X,P(\k)) \times \Omega^l(X,P(\k)) \to \Omega^{k+l}(X), \ \
(a_1,a_2) \mapsto (a_1 \wedge a_2) .$$
Choose a metric on $X$ and let $*$ denote the associated Hodge star
operator $ \Omega^k(X,P(\k)) \to \Omega^{2-k}(X,P(\k)).$ Let 
$$\A(P) = \Omega^1(X,P(\k))$$ 
the affine space of connections and
$$K(P) = \Aut_K(P), \ G(P) = \Aut_G(P(G)) $$
the group of unitary, resp. complex gauge transformations.  For any $A
\in \A(P)$, let $ F_A \in \Omega^2(X,P(\k)) $ denote its curvature.
Yang-Mills theory is the area-dependent quantum field theory with
partition function given 
$$ Z(X) = \sum_{P} Z(P) $$
where the sum is over isomorphism classes of principal $K$-bundles $P$
and $Z(P)$ is defined formally by the path integral
$$ \text{``} Z(P) = \frac{1}{\Vol(K(P))} \int_{\A(P)} \exp( - 
S(A) ) DA \text{''}, \ \ \ S(A) = \frac{1}{2 \eps} \int_X 
(F_A \wedge * F_A) .$$
Formally $Z(P)$ is the pairing of the Duistermaat-Heckman measure for
the action of $K(P)$ on $\A(P)$ with a Gaussian on
$\Omega^2(X,P(\k))$.

A definition of the two-dimensional Yang-Mills integral, including
observables, is given by Levy \cite{le:ym}.  Levy's approach is to
embed the space of connections mod gauge equivalence into the space of
maps of the loop space on $\Sigma$ to the group mod conjugacy via the
holonomy map.  Levy constructs a probability measure on this
``thickening'' of the space of connections and proves that the
Yang-Mills integral is given by the Migdal formula.

Here we will define the Yang-Mills integral by assuming that
localization for the norm-square \eqref{Hform} holds.  The strategy of
defining path integrals by expanding over critical components of the
integrand appears in many places, such as perturbative Chern-Simons
theory \cite{as:cs2}.  The purpose of this section is to show that
with this definition, the Yang-Mills integral is given by the Migdal
formula, and hence agrees with Levy's definition. This might be seen
as an easy two-dimensional analog of the much harder conjecture
regarding the three-dimensional Chern-Simons path integral, that the
``exact'' definition via Reshetikhin-Turaev agrees with the
``perturbative'' definition of Axelrod-Singer.

The action of $K(P)$ on $\A(P)$ is Hamiltonian with moment map minus
the curvature, and so the Yang-Mills function $S(A)$ is the
norm-square of the moment map.  The critical points of $S(A)$ are the
connections satisfying the {\em Yang-Mills equation}
$$ d_A^* F_A = 0 .$$
These are the connections with constant central curvature.  Each such
connection is gauge equivalent to a connection $A'$ with $F_{A'} = *
\xi$; let $\Xi(P)$ denote the set of $\xi$.  In the case $K = U(r)$,
$P$ is the principal $U(r)$ bundle of rank $r$ and degree $d$ over a
surface $X$ of genus at least one,
$$ \Xi(P) = \left\{ (\mu_1,\ldots,\mu_1,\mu_2,\ldots,\mu_2,\ldots,
 \mu_r) \right\} \subset \mathbb{Q}^r $$
the set of non-increasing sequences such that $\mu_j = d_j/r_j$ for
some integers $d_j$ and $r_j$ such that
$ \sum_j d_j = d , \ \ \ \sum_j r_j = r $
and each $\mu_j$ appears $r_j$ times.  If $X$ is genus zero then only
integral $\mu_j$ appear.

Minus the gradient flow of $S(A)$ induces a decomposition of $\A(P)$
into stable manifolds
$$ \A(P) = \bigcup_{\xi \in \Xi} \A(P)_\xi $$
By results of Donaldson \cite{do:ne}, Daskalopolous\cite{da:st},
R{\aa}de \cite{ra:th}, and Atiyah-Bott \cite{at:mo} this is identical
to the decomposition by Harder-Narasimhan type of the corresponding
holomorphic $G$-bundle.  For $K = U(n)$, $\A(P)_\xi$ consists of
connections such that the corresponding holomorphic bundle has
Harder-Narasimhan quotients with ranks $r_j$ and degrees $d_j$.

For each $\xi \in \Xi(P)$, the universal quotient of $\A(P)_\xi$ by
$G(P)$ is the moduli space $\M(X,K_\xi;\xi)$ of $K_\xi$-bundles with
constant central curvature $\xi$.  Define a vector bundle $\nu_\xi \to
\M(X,K_\xi;s\xi)$ by
$$ (\nu_\xi)_{[A]} = (H^1 - H^0)(\ol{\partial}_A,\g/\g_\xi) \oplus
\p_\xi/\g_\xi$$
where $\ol{\partial}_A$ is the corresponding Dolbeault operator.
Except for the factor $\p_\xi/\g_\xi$, this is the virtual normal
bundle for the embedding of moduli stacks induced by $G_\xi \to G$,
see \cite{te:in}.  The inclusion of $\p_\xi/\g_\xi$ in the definition
has to do with the fact that the generic complex automorphism group of
a bundle of type $\xi$ is the corresponding parabolic $P_\xi$, which
means that $- \p_\xi/\g_\xi$ appears in the stacky normal bundle
$\nu_\xi$ but not in the corresponding formula in Section 5.

Let $K'_\xi$ denote the identity component of the generic automorphism
group for $M(X,K_\xi;\xi)$, $K_\xi'' = K_\xi/K'_\xi$.  Let
$K_{\xi,M}''$ denote the (finite) subgroup of $K_\xi''$ contained in
the generic automorphism group. Let $ \M(X,x,K_\xi;\xi)$ denote the
moduli space of bundles with framing at a base point $x$. If every
point in $\M(X,K_\xi;\xi)$ has automorphism group $K_\xi'$ then
$M(X,x,K_\xi;\xi)$ is a locally free $K''_\xi$-space with quotient $
\M(X,K_\xi;\xi)$.  For any $h \in \S(\k_\xi)^{K_\xi}$, let
$\kappa_\xi(h) \in \H_{K'_\xi}(\M(X,K_\xi;\xi))$ denote the
corresponding characteristic class.  Let $\L_{(\xi)}$ denote the
$K_\xi'$-equivariant Liouville form on $ \M(X,K_\xi;\xi)$, with
constant moment map with value $\xi$.  Let $\mu_{\A(X),\xi} \in
\S'(\k^*)^K$ denote the distribution defined by
$$ (\mu_{\A(X),\xi},h) = \int_{ \M(X,K_\xi;\xi) \times \k'_\xi}
\L_{(\xi)} \wedge \Eul(\nu_\xi)_\xi^{-1} \wedge
\kappa_\xi(\Vol_{K_\xi}^K \Res_{K_\xi}^K h),
$$
times $\Vol(K''_\xi/K''_{\xi,M})$, compare with \ref{Hsmooth}.  Let
$$ \Xi(X) = \bigcup_P \Xi(P) $$
and define the Yang-Mills partition function by 
$$ Z(X) := \sum_{\xi \in \Xi(X)} (\mu_{\A(X),\xi}, h) $$
where $h \in \S(\k^*)^K$ is the Fourier transform of $\hat{h}(\zeta) =
\exp \left( - \frac{\eps}{2} \Vert \zeta \Vert^2 \right) .$ (There is
a slight inconsistency with the previous formal definition to the
effect of a missing factor of a power of $\eps$.)

Some care is needed for the definition in the presence of reducible
connections.  Let $\M(X,K)_\nu$ denote the moduli space of flat $K$
bundles on the once-punctured surface, with holonomy around the
puncture conjugate to $\exp(\nu)$.  This space admits a holomorphic
description in terms of semistable bundles with a {\em parabolic
reduction} at the puncture, described in Mehta-Seshadri \cite{ms:pb}.
Let $Z(K)$ denote the center of $K$, and $K'' = K/Z(K)$.  The function
$$Z(X,\nu) := \# \Vol(Z(K)) \Vol(\M(X,K)_\nu)$$
is piecewise polynomial for $\nu \in \k''$.  If every point in
$\M(X,K)$ has automorphism group $Z(K)$ then
$$ \mu_{\A(X),0} = \Vol(K \cdot \nu)^{-1} Z(X,\nu) \mu_{\k^{'',*}} $$
for $\nu$ in a neighborhood of $0$. In case $\M(X,K)$ contains
reducibles, this can be taken as the definition of $\mu_0$.  There are
similar definitions for the other distributions $\mu_\xi$ in the
presence of reducible connections.

The main result of this section is
\begin{theorem} (``Migdal formula'', see \cite[2.51]{wi:qg}) 
\label{migdal}
Let $K$ be a compact connected group.  The $2$-dimensional Yang-Mills
partition function is given by
$$Z(X) = \Vol(K)^{2g}
 \sum_{\nu} (\dim V_\nu)^{2-2g} \hat{h}(\nu + \rho) $$
where the sum is over dominant $\nu$ in the weight lattice $\Lambda^*$
plus $\rho$.
\end{theorem}

Here $\rho$ is the half-sum of the positive roots which is a weight if
$\k$ is spinnable.  Before we give the proof, we note the corollary
(as already discussed in \cite{wi:qg})

\begin{corollary}  Suppose that $K$ is semisimple and $g \ge 2$.  The volume of the moduli space $\M(X,K)$ is 
$$ \Vol(\M(X,K)) = \#Z(K) \dim(K)^{2g} \sum_{\nu} (\dim
 V_\nu)^{2-2g} $$
where $Z(K)$ is the center of $K$.
\end{corollary}

\begin{proof}  Take the limit $\eps \to 0$ in Theorem \ref{migdal}.
By definition of $Z(X)$, the limit
$$ \lim_{\eps \to 0} Z(X) = \#Z(K)^{-1} \Vol(\M(X,K)) .$$
On the other hand, the limit of the right hand side of \ref{migdal} is
$$ \dim(K)^{2g} \sum_{\nu} (\dim V_\nu)^{2-2g} ,$$ 
which proves the corollary.
\end{proof}

The measures $\mu_{\A(X),\xi}$ for $\xi$ generic can be described as
follows. The moduli space $\M(X,K_\xi,\xi)$ is the Jacobian of torus
bundles with first Chern class $\xi$, and is diffeomorphic to
$T^{2g}$.  The characteristic classes of the bundle $\nu$ are computed
in \cite{te:kt},\cite[p.8]{te:in}.  One obtains
$$\Eul(\nu_\xi)= (-1)^{2\rho(\xi)} |\Eul(\k/\t)|^{2g-2} .$$
Integrating over $\M(X,K_\xi,\xi)$ gives
\begin{eqnarray*} \label{mut}
\mu_{\A(X),\xi} &=&  i^{(2g-1)\dim(K/T)/2}
(-1)^{2\rho(\xi)}\Ind_T^K \int_{T^{2g}}
\exp(\omega_{(\xi)}) \delta_\xi \Eul(\k/\t)_{\xi}^{1 - 2g}
 \\ &=&  i^{(2g-1)\dim(K/T)/2} (-1)^{2\rho(\xi)} \Ind_T^K \Vol(T^{2g})
\delta(\xi) \Eul(\k/\t)_{\xi}^{1 - 2g}
.\end{eqnarray*}

The proof of Theorem \ref{migdal} is based on the idea, introduced by
C. Teleman \cite{te:kt}, that the sum over strata is the same as the
sum of contributions from the $T$-bundles.  Define
$$ \mu_{\A(X)} := \sum_{\xi \in \Xi(X)} \mu_{\A(X),\xi} \in
\D'(\k^*)^K $$
which is a kind of Duistermaat-Heckman measure for $\A(X)$.  Let
$$\nu_{\A(X),\xi} := \frac{i^{(g-\hh)\dim(K/T)}}{\# W_\xi}
(-1)^{2\rho(\xi)} \Ind_T^K \Vol(T^{2g}) \delta(\xi)
\Eul(\k/\t)_{\zeta}^{1 - 2g} $$
if $\M(X,K_\xi,\xi)$ contains $T$-bundles, and zero otherwise.  Here
  $\zeta \in \t_+^*$ is any regular element.  We wish to compare
$\mu_{\A(X)}$ with 
$$ \nu_{\A(X)} := \sum_{\xi \in \Xi(X)} \nu_{\A(X),\xi} $$
which is the sum of the ``fixed-point contributions'' from
$T$-bundles.  For any distribution $\mu \in \D'(\k)^K$, define a
distribution $\mu_T \in \D'(\t)^{\on{sign}(W)}$ by
$$ (\mu_T, \Vol_T^K \Res_T^K h ) := (\mu,h) .$$
The map $\mu \mapsto \mu_T$ is a right inverse to $\Ind_T^K$.  We will
need the following lemma:

\begin{lemma} \label{supplem}
\begin{enumerate}
\item $ \mu_{\A(X),T}$ is invariant under translation by the coweight
lattice $\Lambda$ and anti-invariant under $W$.  (In other words,
anti-invariant under the action of the affine Weyl group.)
\item $ (\mu_{\A(X),\xi} - \nu_{\A(X),\xi})_T$ has Fourier transform
supported in $\t_{\sing}$.
\end{enumerate}
\end{lemma} 

The lemma is obtained by taking the high level limit of the
corresponding $K$-theoretic statements in \cite{te:in}.  Since
$\mu_{\A(X),T}$ is a periodic distribution, its Fourier transform
$\F_\t^{-1} \mu_{\A(X),T}$ is a sum of delta functions at weights:
$$ \F_\t^{-1} \mu_{\A(X),T} = \sum c_\lambda \delta_\lambda .$$
Since $\mu_{\A(X),T}$ is $W$-anti-invariant, $c_\lambda = 0 $ unless
$\lambda$ is regular.  By part (b) of Lemma \ref{supplem}, $
\mu_{\A(X),T}$ is equal to $\nu_{\A(X),T} $ plus a distribution
whose Fourier transform is supported in $\t_{\sing}$.  We have
$$ \nu_{\A(X),T} = \frac{i^{(g-\hh)\dim(K/T)}}{\# W} \sum_{\xi
 \in \Lambda} \delta_\xi \Vol(T^{2g}) (-1)^{2\rho(\xi)} \Eul(\k/\t)^{1
 - 2g} .$$
By the Poisson summation formula
$$ \F_t^{-1} \mu_{\A(X),T} = i^{-\dim(K/T)/2} \sum_{\lambda \in
\Lambda^* + \rho} (\# W)^{-1} \delta_{\lambda} \Vol(T^{2g})
\prod_{\alpha > 0} 2\pi (\alpha,\lambda)^{1 - 2g} $$
where (since $c_\lambda = 0$ for singular $\lambda$) the sum is over
{\em regular} $\lambda$.  Hence
$$ \mu_{\A(X)} = i^{-\dim(K/T)/2} \Vol(K)^{2g-1} \Vol(T) \Ind_T^K
\sum_{\lambda} \delta_{\lambda + \rho} \dim(V_\lambda)^{1-2g} $$
where the sum is over $\lambda$ such that $\lambda + \rho$ is a
dominant weight.  (If $\lambda$ is not a weight, $V_\lambda$ is a
representation of the universal cover of $K$.)  Finally pairing with
the Gaussian $h$ gives
\bea ( \mu_{\A(X)}, h)
&=& 
\frac{\Vol(K)^{2g-1} \Vol(T)}{
 i^{\dim(K/T)/2} }
 \left( \sum_{\lambda} \delta_{\lambda +
  \rho}\dim(V_\lambda)^{1-2g} , \Vol_T^K \Res_T^K \hat{h} \right) \\
&=& \left(  \Vol(K)^{2g} \sum_{\lambda}  \delta_{\lambda + \rho} \dim(V_\lambda)^{2-2g} ,
\Res_T^K \hat{h} \right) \\ 
&=& \Vol(K)^{2g} \sum_{\lambda} \dim(V_\lambda)^{2-2g} \hat{h}(\lambda +
 \rho) \eea
which completes the proof of \ref{migdal}.  This computation is done
on a physics level of rigor by Blau-Thompson \cite{bl:lo}.  The main
point is that the contribution of the semistable stratum is not affine
Weyl-invariant, but only becomes so after adding the contributions
from the higher strata.  As in \cite{te:in}, the additional symmetry
removes the necessity of doing any hard computations, that is, any
integrals other than integrals over Jacobians.

\begin{example}  Let $K = SU(2)$ and identify $\t \to \R$
so that the weight lattice is $\Z/2$ and coweight lattice $\Z$.  If
$X$ has genus $g = 1$, an explicit computation shows
$$  Z(X,\nu) =  \qq \Vol(T^2) (1 -  2 \nu) .$$
For $\xi$ a positive integer, $K_\xi = K$ and $M(X,T;\xi) = T^{2}$.
The normal bundle $\nu_\xi = (\k/\t)^2$, hence
\begin{eqnarray*} \mu_{\A(X),\xi} &=& 
\hh \Ind_T^K \Vol(T^{2}) \delta(\xi) \Eul(\k/\t)_+^{-1} -
\delta(-\xi) \Eul(\k/\t))_-^{-1} \\ 
&=& \hh \Vol(T^{2}) ( \delta(\xi) H_+ -
\delta(-x) H_-) \end{eqnarray*}
where $H_\pm$ are the Heaviside distributions.  The sum is the
sawtooth distribution shown below in solid lines in Figure \ref{Z}.
The dotted line is the contribution from $\xi = 0$.  
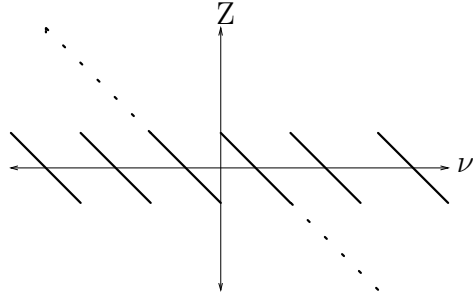
\begin{figure}[h]
\setlength{\unitlength}{0.00030333in}
\begingroup\makeatletter\ifx\SetFigFont\undefined%
\gdef\SetFigFont#1#2#3#4#5{%
  \reset@font\fontsize{#1}{#2pt}%
  \fontfamily{#3}\fontseries{#4}\fontshape{#5}%
  \selectfont}%
\fi\endgroup%
{\renewcommand{\dashlinestretch}{30}
\begin{picture}(7902,4747)(0,-10)
\path(3652.000,4402.000)(3622.000,4522.000)(3592.000,4402.000)
\path(3622,4522)(3622,22)
\path(3592.000,142.000)(3622.000,22.000)(3652.000,142.000)
\thicklines
\path(3622,2722)(4822,1522)
\path(4822,2722)(6022,1522)
\path(6322,2722)(7522,1522)
\path(2422,2722)(3622,1522)
\path(1222,2722)(2422,1522)
\path(22,2722)(1222,1522)
\thinlines
\path(142.000,2152.000)(22.000,2122.000)(142.000,2092.000)
\path(22,2122)(7222,2122)(7522,2122)
\path(7402.000,2092.000)(7522.000,2122.000)(7402.000,2152.000)
\thicklines
\dashline{90.000}(4822,1522)(6322,22)
\dashline{90.000}(2422,2722)(622,4522)(622,4447)
\put(7672,2047){\makebox(0,0)[lb]{\smash{{{\SetFigFont{12}{14.4}{\rmdefault}{\mddefault}{\updefault}$\nu$}}}}}
\put(3547,4597){\makebox(0,0)[lb]{\smash{{{\SetFigFont{12}{14.4}{\rmdefault}{\mddefault}{\updefault}Z}}}}}
\end{picture}
}
\caption{One-point function for genus one, $K = SU(2)$ \label{Z}}
\end{figure}

\end{example} 

It seems an interesting question whether a similar definition could be
used for other path integrals, for instance, holomorphic Yang-Mills
theory in four dimensions.  On the other hand, other path integrals
such as full four-dimensional Yang-Mills or two-dimensional Yang-Mills
with observables do not seem to admit heuristic interpretations as
pairings in equivariant cohomology, and it appears unlikely that the
techniques described here would apply.

\appendix

\newcommand{\tq}{{\f{3}{4}}}

\section{Convergence of the gradient flow}

This section contains a proof that the norm-square of minus the
gradient flow of the moment map converges.  This result appeared some
time ago in an unpublished manuscript by Duistermaat, who used the
gradient inequality of Lojasiewicz \cite{lo:en}, see Lerman
\cite{le:gra}.  Here we prove the necessary gradient inequality
directly, using the local model \eqref{model}, and obtain explicit
estimates for the rate of convergence.  The same estimates were
established for the infinite dimensional Yang-Mills heat flow by
R{\aa}de \cite{ra:ym}; our goal is to put the finite dimensional case
on equal footing.

First we discuss some background on gradient flows.  Let $M$ be a
Riemannian manifold with metric $g$, $f \in C^\infty(M)$ and
$\grad(f)$ its gradient.  Let $\crit(f)$ the critical set of $f$.  For
any $x \in \crit(f)$, let $L(f,x)$ be the set of $\gamma > 0$ such
that there exists a neighborhood $U$ of $x$ and a constant $C > 0$ so
that for all $m \in U$,
$$ \Vert \grad(f)(m) \Vert > C|f(m) - f(x)|^{\gamma} .$$

\begin{theorem}[Lojasiewicz gradient inequality] \label{loj}
If $f:\R^n \to \R$ is a real-analytic function then for every critical
point $x$ of $f$ there exists a $\gamma \in L(f,x)$ with $\gamma \in
[\hh,1)$.
\end{theorem} 
\noindent The Riemannian metric on $\R^n$ is assumed to be the
standard one.  However, the same inequality holds for an arbitrary
metric, with the same exponent but possibly different constant.
Therefore, the same result holds for functions on Riemannian manifolds
that are real-analytic near their critical sets.  We will not use
Theorem \ref{loj} in this paper but instead a special case which is
easy to prove:
\begin{lemma} \label{homog} Let 
$f:\R^n \to \R$ be a homogeneous polynomial of degree $d$.  $ 1- 1/d
\in L(f,0)$.
\end{lemma}
\begin{proof}  Let $v \in \R^n$ have norm $1$ and 
$f_v(t) = f(tv)$.  Then 
$$\Vert \grad(f)(tv) \Vert \ge | f_v'(t) | = (d-1) | f_v(t) |^{1 - 1/d} = (d-1) | f(tv) |^{1 - 1/d} .$$
\end{proof} 

The following theorems are probably well-known.  Suppose that $f$ is
proper and bounded from below.  Since $f^{-1}(-\infty,c]$ is compact
for any $c > 0$, the flow $\varphi_t:M \to M$ of $-\grad(f)$ is
defined for all times $t$.
\begin{theorem}  \label{grad}  
Let $c$ be a critical value of $f$.  Suppose that there exists $\gamma
\in (0,1)$ such that $\gamma \in L(f,x)$ for every $x \in \crit(f)
\cap f^{-1}(c)$.
\begin{enumerate}
\item Any trajectory $m_t$ of $-\grad(f)$ such that $f(m_t) \to c$ has
a unique limit $m_\infty$ as $t \to \infty$.
\item Let $W_c$ denote the {\em stable set} of points $m \in M$ with
$m_\infty \in f^{-1}(c)$.  The map $m \to m_\infty$ is a deformation
retraction of $W_c$ onto $f^{-1}(c) \cap \crit(f)$.
\end{enumerate}
\end{theorem}

\begin{theorem} \label{converge}
For any $m \in M$ and $\gamma \in L(f,m_\infty)$, there exist
constants $C,k$ and a time $T$ such that if $t > T$ then
\begin{enumerate}
\item if $\gamma \in (\hh,1)$ then $d(m_t,m_\infty) \leq Ct^{(\gamma
  - 1)/(2\gamma - 1)},$ and
\item if $\gamma = \hh$ then  $d(m_t,m_\infty) \leq Ce^{-kt}$.
\end{enumerate}
\end{theorem} 
The following is a sketch of proof, see also \cite{le:gra}.  Suppose
that $\gamma \in (0,1)$ lies in $L(f,x)$ for all $x \in \crit(f)$.
Since $f$ is proper, for each critical value $c$ of $f$ there exists a
neighborhood $U_c$ of $f^{-1}(c) \cap \crit(f)$ such that if $m \in
U_c$ then
\begin{equation} \label{gradineq}
 \Vert \grad(f)(m) \Vert \ge ( f(m) - c)^\gamma .\end{equation}
Let $m \in M$ and $m_t$ the trajectory of $-\grad(f)$ with $m(0) =
m$.  Since $f(m_t)$ is non-increasing, all limit points of $m_t$ lie
in $\crit(f) \cap f^{-1}(c)$ for some $c$.  Using properness of $f$
again, there exists $T > 0$ such that $m_t$ lies in some $U_c$ for $t
> T$.  For $t > T$, 
\begin{eqnarray} \label{ddt}
 \ddt (f(m_t) - c )^{1 - \gamma} &=& - \Vert \grad(f)(m_t) \Vert^2 (
 f(m_t) - c )^{-\gamma} \\ \label{ddt2}
&\geq& - C \Vert \grad(f)(m_t) \Vert
 \end{eqnarray}
by \eqref{gradineq}.  Hence for $t_1,t_2 > T$,
\begin{eqnarray*} d(m(t_1),m(t_2)) &\leq& 
\int_{t_1}^{t_2} \d t \Vert \grad(f(m_t) \Vert\\
&\leq& C ( ( f(m(t_1) - c )^{1 - \gamma} - ( f(m(t_2)) - c)^{1 -
  \gamma}) .\end{eqnarray*}
By the Cauchy criterion, $m_t$ converges to a critical point $m_\infty
\in f^{-1}(c)$.  \eqref{ddt},\eqref{ddt2} also imply that for $t > T$
\begin{equation} \label{further}
 d(m_t,m_\infty) \leq C (f(m_t) - c )^{1 - \gamma} .\end{equation}
Next we show that the map $W_c \to f^{-1}(c) \cap \crit(f), \ m
\mapsto m_\infty $ is a deformation retraction, that is, that $ W_c
\times [0,\infty] \to W_c, \ m \mapsto m_t $ is continuous.  Let $\eps
> 0$.  By \eqref{further}, there exists $\delta > 0 $ so that $f(m) -
c < \delta$ and $m \in W_c$ imply $d(m,m_\infty) < \eps/3$.  Fix $m_1
\in W_c$ and let $t$ be sufficiently large so that $f(m_{1,t}) - c <
\delta$.  By smooth dependence on initial conditions, if $m_2$ is
sufficiently close to $m_1$ then $d(m_{1,t},m_{2,t}) < \eps/3$.  For
$\eps$ sufficiently small, $f(m_{2,t}) - c< \delta$.  Hence if $m_2$
lies in $W_c$ and is sufficiently close to $m_1$ then
\begin{eqnarray*} 
d(m_{1,\infty},m_{2,\infty}) &\leq& d(m_{1,\infty},m_{1,t}) +
d(m_{1,t},m_{2,t}) + d(m_{2,t},m_{2,\infty}) \\ &<& \eps/3 + \eps/3 +
\eps/3 = \eps  \end{eqnarray*}
and also $d(m_{2,t},m_{1,\infty}) < \eps$.  This completes the proof
of Theorem \ref{grad}.

Going back to \eqref{gradineq} we have 
$$ \ddt (f(m_t) - c) = - \Vert \grad(f)(m) \Vert^2 \leq -C (f(m_t) -
c)^{2 \gamma} .$$
Integrating gives
$$ f(m_t) - c \leq \left\{ 
\begin{array}{ll} C(t-T)^{-1/(2\gamma - 1)} & \gamma \in (\hh,1) \\
                  Ce^{-k(t - T)} & \gamma = \hh \end{array} \right\}.  $$ 
Using \eqref{further} this proves Theorem \ref{converge}

We apply Theorems \ref{grad}, \ref{converge} to the norm-square of the
moment map.  Let $K$ be a compact Lie group with Lie algebra $\k$, and
$(\ , \ ): \k \times \k \to \R$ an invariant inner product.  From now
on, we use the inner product to identify $\k$ with its dual $\k^*$.
Let $M$ be a Hamiltonian $K$-manifold with proper moment map $\Phi: M
\to \k$ and $f = \hh(\Phi,\Phi) .$ Choose an invariant compatible
almost complex structure $J$, so that $\omega(\cdot,J \cdot)$ defines
a $K$-invariant Riemannian metric $g$ on $M$. Let
$$M = \bigcup_{\xi \in \Xi(M)} M_\xi$$ 
denote the Kirwan-Ness stratification into stable sets for
$-\grad(f)$.  For each $\xi \in \Xi(M)$, let $C_\xi$ denote the
corresponding component (not necessarily connected) of $\crit(f)$, and
$Z_\xi$ the union of components of the fixed point set of $U(1)_\xi$
meeting $\Phinv(\xi)$.

\begin{lemma} \label{lb} Let $x$ be a point in
$\Phinv(\xi) \cap \crit(f)$.  (a) $\tq \in L(f,x)$.  (b) If $x$ lies in
the principal orbit-type stratum of $Z_\xi$ then $\hh \in L(f,x)$.
\end{lemma}

We apply the local model \eqref{model} to $y = x \in \crit(f)$, so
that $\xi \in \k_x$.  We identify a neighborhood of $x$ in $M$ with a
neighborhood of $0$ in $\k_x^\perp \oplus T_0 S \oplus (\k_\xi^* \cap
\k_x^\circ)$.  We denote by $g_0$ the constant metric that is the
direct sum of the metric on $T_0 S$ and the metrics on $\k_x^\perp
\oplus (\k_\xi^* \cap \k_x^\circ)$ induced by the metric on $\k$. Let
$\Vert \cdot \Vert_0$ denote the corresponding norm; we will denote by
the same notation its restriction to $T_0 S$.  Let $m \in M$ and $[
k,s,\nu] $ be the image of $m$ in $U'$.  Since $f$ is $K$-invariant,
it suffices to consider the case $k = 1$.  For $m$ sufficiently close
to $x$ we have
\begin{eqnarray*} \Vert \grad(f)(m) \Vert &=& 
\Vert \Phi(m)_M(m) \Vert \\
&\ge& C_1 \Vert \Phi(m)_M(m) \Vert_0 \\
&=&   C_1 \Vert (\nu, (\Phi_S(s) + \xi)\cdot s,
(\Phi_S(s) + \xi) \cdot \nu) \Vert_0 \\
&=&  C_1 ( \Vert (\xi + \Phi_S(s)) \cdot s \Vert_0^2 + \Vert \nu
     \Vert^2)^{1/2} . \end{eqnarray*}
Let $S_0 \subset S$ denote the fixed point set of $K_x$, $S_1$ the
symplectic complement of $S_0$ in the fixed point set of $K_\xi$, and
$S_2$ the symplectic complement of $S_1 \oplus S_2$ so that $ S = S_0
\oplus S_1 \oplus S_2 .$ Since $\Phi$ vanishes on $S_0$, $f$ and
$\grad(f)$ are independent of $s_0 \in S_0$.  $f$ is homogeneous of
degree $4$ on $S_1$.  By Lemma \ref{homog}
\begin{eqnarray*}
 \Vert \Phi_S(s_0,s_1,0) \cdot (s_0,s_1,0) \Vert_0^2 &=&
\Vert \grad (\Phi_S,\Phi_S)(s_0,s_1,0)/2 \Vert_0^2 \\
&\ge& C_2 \Vert \Phi_S(s_0,s_1,0) \Vert^{3}
.\end{eqnarray*}
We expand
\begin{multline*}
\Vert (\xi + \Phi_S(s)) \cdot s \Vert_0^2 = \Vert \Phi_S(s_0,s_1,0)
\cdot (s_0,s_1,0) \Vert_0^2 +\\
 2g_0( \Phi_S(0,0,s_2) \cdot (s_0,s_1,0),
\Phi_S(s_0,s_1,0) \cdot (s_0,s_1,0))+  \Vert (\xi + \Phi_S(s_0,s_1,0))
\cdot s_2 \Vert_0^2  \\
+  \Vert \Phi_S(0,0,s_2) \cdot (s_0,s_1,0) \Vert_0^2 + 
2 g_0( (\xi + \Phi_S(s_0,s_1,0))  \cdot (0,0,s_2) , 
         \Phi_S(0,0,s_2) \cdot (0,0,s_2)) \\ + \Vert
\Phi_S(0,0,s_2) \cdot (0,0,s_2) \Vert_0^2 \end{multline*}
which are terms of degree $0,2,2,4,4,6$ respectively in $s_2$.
Because $\xi$ acts on $S_2$ without fixed points, the terms of degree
$2$ sum to a positive quadratic form for $s_1$ sufficiently small.  It
follows that
$$ \Vert (\xi + \Phi_S(s)) \cdot s \Vert_0^2 \ge C_3 ( \Vert
\Phi_S(s_0,s_1,0) \Vert^{3} + \Vert s_2 \Vert_0^2) $$
for $s_1,s_2$ sufficiently small.  By \eqref{momentf}
\begin{eqnarray*} \label{diff}
 f(m) - f(x) &=& \hh \Vert \Phi_S(s) + \xi \Vert^2 + \hh \Vert
\nu \Vert^2 - \hh \Vert \xi \Vert^2 \\
&=& \hh \Vert \Phi_S(s) \Vert^2 + (\Phi_S(s),\xi) + \hh \Vert \nu
\Vert^2 \\
&=& \hh \Vert \Phi_S(s_0,s_1,0) \Vert^2 +
(\Phi_S(s_0,s_1,0),\Phi_S(0,0,s_2)) \\
&& + \hh \Vert \Phi_S(0,0,s_2)  \Vert^2 + (\xi,\Phi_S(0,0,s_2)) + 
\hh \Vert \nu \Vert^2 \end{eqnarray*}
which is bounded by
$  \hh \Vert \Phi_S(s_0,s_1,0) \Vert^2 + C_4 \Vert s_2 \Vert^2 
+ \hh \Vert \nu \Vert^2 .$
For $a,b \in [0,1] $ we have
%
$ (a^3 + b^2)^2 \ge 
(a^2 + b^2)^3/4. $
%
Applying this inequality with $a = \Vert \Phi_S(s_0,s_1,0) \Vert, \ b
= (\Vert s_2 \Vert_0^2 + \Vert \nu \Vert^2)^{1/2}$ gives
\begin{eqnarray*}
\Vert \grad(f) \Vert & \ge &   
C_5 \left( \Vert \Phi_S(s_0,s_1,0) \Vert^{3} + \Vert s_2 \Vert_0^{2} + \Vert
\nu \Vert^{2} \right)^{\hh} \\
&\geq& C_6 \left( \Vert  \Phi_S(s_0,s_1,0) \Vert^{2} +
\Vert s_2 \Vert_0^{2} + \Vert \nu \Vert^{2} \right)^{\tq} \\
&\geq& C_7 (f(m) - f(x))^{\tq} 
\end{eqnarray*}
which completes the proof of part (a) of the Lemma.  To prove (b),
note that if $x$ is contained in the principal orbit-type stratum of
$Z_\xi$ then $S_1$ is trivial.  Hence for $t > T$
%
$$ \Vert \grad(f) \Vert  \ge    
C \left(\Vert s_2 \Vert_0^{2} + \Vert
\nu \Vert^{2} \right)^{\hh} 
\geq C(f(m_t) - c)^{\hh}. 
$$
%
It now follows from Theorem \ref{grad} 
\begin{theorem} \label{norm}
\begin{enumerate}
\item Every trajectory $m_t$ of $-\grad(f)$ converges to a point
$m_\infty \in \crit(f)$.
\item For all $ \xi \in \Xi(M)$, the map $m \mapsto m_\infty$ is a
deformation retraction of $M_\xi$ onto $C_\xi$.  
\item If $m_\infty$ is contained in the principal orbit-type stratum
of $Z_\xi$ then there exist constants $C,k,T >0$ such that
$d(m_t,m_\infty) < Ce^{-kt}$ for $t > T$.  Otherwise, there exist
constants $C,T>0$ such that $d(m_t,m_\infty) < Ct^{-\hh}$ for 
$t > T$.
\end{enumerate}
\end{theorem}


\def\cprime{$'$} \def\cprime{$'$} \def\cprime{$'$}
  \def\polhk#1{\setbox0=\hbox{#1}{\ooalign{\hidewidth
  \lower1.5ex\hbox{`}\hidewidth\crcr\unhbox0}}}


\begin{thebibliography}{10}

\bibitem{at:mo}
M.~F. Atiyah and R.~Bott.
\newblock The {Y}ang-{M}ills equations over {R}iemann surfaces.
\newblock {\em Phil. Trans. Roy. Soc. London Ser. A}, 308:523--615, 1982.

\bibitem{at:mom}
M.~F. Atiyah and R.~Bott.
\newblock The moment map and equivariant cohomology.
\newblock {\em Topology}, 23(1):1--28, 1984.

\bibitem{as:cs2}
S.~Axelrod and I.~M. Singer.
\newblock Chern-{S}imons perturbation theory. {II}.
\newblock {\em J. Differential Geom.}, 39(1):173--213, 1994.

\bibitem{be:he}
N.~Berline, E.~Getzler, and M.~Vergne.
\newblock {\em Heat Kernels and {D}irac Operators}, volume 298 of {\em
  Grundlehren der mathematischen Wissenschaften}.
\newblock Springer-Verlag, Berlin-Heidelberg-New York, 1992.

\bibitem{be:ze}
N.~Berline and M.~Vergne.
\newblock Z\'ero d'un champ de vecteurs et classes caract\'eristiques
  \'equivariantes.
\newblock {\em Duke Math. J.}, 50:539--549, 1983.

\bibitem{bi:so}
A.~Bia{\l}ynicki-Birula.
\newblock Some theorems on actions of algebraic groups.
\newblock {\em Ann. of Math. (2)}, 98:480--497, 1973.

\bibitem{bl:lo}
M.~Blau and G.~Thompson.
\newblock Localization and diagonalization: a review of functional integral
  techniques for low-dimensional gauge theories and topological field theories.
\newblock In {\em Functional integration (Carg\`ese, 1996)}, volume 361 of {\em
  NATO Adv. Sci. Inst. Ser. B Phys.}, pages 363--410. Plenum, New York, 1997.

\bibitem{ca:ko}
A.~Canas~da Silva and V.~Guillemin.
\newblock On the {K}ostant multiplicity formula for group actions with
  non-isolated fixed points.
\newblock {\em Adv. Math.}, 123(1):1--15, 1996.

\bibitem{da:st}
G.~D. Daskalopoulos.
\newblock The topology of the space of stable bundles on a compact {R}iemann
  surface.
\newblock {\em J. Differential Geom.}, 36(3):699--746, 1992.

\bibitem{do:ne}
S.~K. Donaldson.
\newblock A new proof of a theorem of {N}arasimhan and {S}eshadri.
\newblock {\em J. Differential Geom.}, 18(2):269--277, 1983.

\bibitem{du:eq}
J.~J. Duistermaat.
\newblock Equivariant cohomology and stationary phase.
\newblock In {\em Symplectic geometry and quantization, (Sanda and Yokohama,
  1993)}, volume 179 of {\em Contemp. Math.}, pages 45--62, Providence, RI,
  1994. Amer. Math. Soc.

\bibitem{freed:comp}
D.~S. Freed and R.~E. Gompf.
\newblock Computer calculation of {W}itten's {$3$}-manifold invariant.
\newblock {\em Comm. Math. Phys.}, 141(1):79--117, 1991.

\bibitem{gu:on}
V.~Guillemin, E.~Lerman, and S.~Sternberg.
\newblock On the {K}ostant multiplicity formula.
\newblock {\em J. Geom. Phys.}, 5(4):721--750, 1988.

\bibitem{gu:symu}
V.~Guillemin, E.~Lerman, and S.~Sternberg.
\newblock {\em Symplectic fibrations and multiplicity diagrams}.
\newblock Cambridge University Press, Cambridge, to appear.

\bibitem{gu:no}
V.~Guillemin and S.~Sternberg.
\newblock A normal form for the moment map.
\newblock In S.~Sternberg, editor, {\em Differential Geometric Methods in
  Mathematical Physics}, volume~6 of {\em Mathematical Physics Studies}, pages
  161--175, Jerusalem, 1982, 1984. D. Reidel Publishing Company, Dordrecht.

\bibitem{gu:eqdr}
V.~W. Guillemin and S.~Sternberg.
\newblock {\em Supersymmetry and equivariant de {R}ham theory}.
\newblock Springer-Verlag, Berlin, 1999.
\newblock With an appendix containing two reprints by Henri Cartan [MR {\bf
  13},107e; MR {\bf 13},107f].

\bibitem{ha:rd}
R.~Hartshorne.
\newblock {\em Residues and duality}.
\newblock Lecture notes of a seminar on the work of A. Grothendieck, given at
  Harvard 1963/64. With an appendix by P. Deligne. Lecture Notes in
  Mathematics, No. 20. Springer-Verlag, Berlin, 1966.

\bibitem{he:re}
P.~Heinzner and F.~Loose.
\newblock Reduction of complex {H}amiltonian ${G}$-spaces.
\newblock {\em Geom. Funct. Anal.}, 4(3):288--297, 1994.

\bibitem{ho:an2}
L.~H{\"o}rmander.
\newblock {\em The analysis of linear partial differential operators. {I}{I}}.
\newblock Springer-Verlag, Berlin, 1983.
\newblock Differential operators with constant coefficients.

\bibitem{je:lo1}
L.~C. Jeffrey and F.~C. Kirwan.
\newblock Localization for nonabelian group actions.
\newblock {\em Topology}, 34:291--327, 1995.

\bibitem{je:in}
L.~C. Jeffrey and F.~C. Kirwan.
\newblock Intersection theory on moduli spaces of holomorphic bundles of
  arbitrary rank on a {R}iemann surface.
\newblock {\em Ann. of Math. (2)}, 148(1):109--196, 1998.

\bibitem{ki:coh}
F.~C. Kirwan.
\newblock {\em Cohomology of Quotients in Symplectic and Algebraic Geometry},
  volume~31 of {\em Mathematical Notes}.
\newblock Princeton Univ. Press, Princeton, 1984.

\bibitem{ki:con}
F.~C. Kirwan.
\newblock Convexity properties of the moment mapping, {III}.
\newblock {\em Invent. Math.}, 77:547--552, 1984.

\bibitem{ku:eq}
S.~Kumar and M.~Vergne.
\newblock Equivariant cohomology with generalized coefficients.
\newblock {\em Ast{\'e}risque}, 215:109--204, 1993.

\bibitem{le:gra}
E.~Lerman.
\newblock Gradient flow of the norm squared of a moment map.
\newblock math.SG/0410568.

\bibitem{le:co}
E.~Lerman, E.~Meinrenken, S.~Tolman, and C.~Woodward.
\newblock Non-abelian convexity by symplectic cuts.
\newblock {\em Topology}, 37:245--259, 1998.

\bibitem{le:ym}
T.~L{\'e}vy.
\newblock Yang-{M}ills measure on compact surfaces.
\newblock {\em Mem. Amer. Math. Soc.}, 166(790):xiv+122, 2003.

\bibitem{lo:en}
S.~{L}ojasiewicz.
\newblock Sur les ensembles semi-analytiques.
\newblock In {\em Actes du Congr\`es International des Math\'ematiciens (Nice,
  1970), Tome 2}, pages 237--241. Gauthier-Villars, Paris, 1971.

\bibitem{ma:vo}
C.-M. Marle.
\newblock Le voisinage d'une orbite d'une action hamiltonienne d'un groupe de
  {L}ie.
\newblock In {\em South Rhone seminar on geometry, II (Lyon, 1983)}, Travaux en
  Cours, pages 19--35. Hermann, Paris, 1984.

\bibitem{mar:sy}
S.~Martin.
\newblock Symplectic quotients by a nonabelian group and by its maximal torus.
\newblock math.SG/0001002.

\bibitem{ms:pb}
V.~B. Mehta and C.~S. Seshadri.
\newblock Moduli of vector bundles on curves with parabolic structure.
\newblock {\em Math. Ann.}, 248:205--239, 1980.

\bibitem{me:onwi}
E.~Meinrenken.
\newblock On {W}itten's formulas for intersection pairings on the moduli space
  of flat {G}-bundles.
\newblock 2003 preprint, in preparation.

\bibitem{me:sym}
E.~Meinrenken.
\newblock Symplectic surgery and the {S}pin$^{\rm c}$-{D}irac operator.
\newblock {\em Adv. in Math.}, 134:240--277, 1998.

\bibitem{ne:st}
L.~Ness.
\newblock A stratification of the null cone via the moment map.
\newblock {\em Amer. J. Math.}, 106(6):1281--1329, 1984.
\newblock with an appendix by D. Mumford.

\bibitem{pa:loc}
P.-E. Paradan.
\newblock Formules de localisation en cohomologie equivariante.
\newblock {\em Compositio Math.}, 117(3):243--293, 1999.

\bibitem{pa:mo}
P.-E. Paradan.
\newblock The moment map and equivariant cohomology with generalized
  coefficients.
\newblock {\em Topology}, 39(2):401--444, 2000.

\bibitem{pa:lo}
P.-E. Paradan.
\newblock Localization of the {R}iemann-{R}och character.
\newblock {\em J. Funct. Anal.}, 187(2):442--509, 2001.

\bibitem{pr:dh}
E.~Prato and S.~Wu.
\newblock Duistermaat-{H}eckman measures in a non-compact setting.
\newblock {\em Compositio Math.}, 94(2):113--128, 1994.

\bibitem{ra:ym}
J.~R{\aa}de.
\newblock On the {Y}ang-{M}ills heat equation in two and three dimensions.
\newblock {\em J. Reine Angew. Math.}, 431:123--163, 1992.

\bibitem{ra:th}
A.~Ramanathan.
\newblock Moduli for principal bundles over algebraic curves. {I}.
\newblock {\em Proc. Indian Acad. Sci. Math. Sci.}, 106(3):301--328, 1996.

\bibitem{rt:in}
N.~Reshetikhin and V.~G. Turaev.
\newblock Invariants of {$3$}-manifolds via link polynomials and quantum
  groups.
\newblock {\em Invent. Math.}, 103(3):547--597, 1991.

\bibitem{sh:gl}
M.~Shub, A.~Fathi, and R.~Langevin.
\newblock {\em Global stability of dynamical systems}.
\newblock Springer-Verlag, New York, 1987.
\newblock Translated from the French by Joseph Christy.

\bibitem{te:kt}
C.~Teleman.
\newblock {K}-theory of the moduli of bundles over a {R}iemann surface and
  deformations of the {V}erlinde algebra.
\newblock math.AG/0306347, to appear in the {S}egalfest proceedings.

\bibitem{te:qu}
C.~Teleman.
\newblock The quantization conjecture revisited.
\newblock {\em Ann. of Math. (2)}, 152(1):1--43, 2000.

\bibitem{te:in}
C.~Teleman and C.~T. Woodward.
\newblock The index formula on the moduli of ${G}$-bundles.
\newblock math.AG/0312154, submitted.

\bibitem{ve:mu}
M.~Vergne.
\newblock Multiplicities formula for geometric quantization {I}, {II}.
\newblock {\em Duke Math. J.}, 82:143--194, 1996.

\bibitem{wi:qg}
E.~Witten.
\newblock On quantum gauge theories in two dimensions.
\newblock {\em Comm. Math. Phys.}, 141:153--209, 1991.

\bibitem{wi:tw}
E.~Witten.
\newblock Two-dimensional gauge theories revisited.
\newblock {\em J. Geom. Phys.}, 9:303--368, 1992.

\end{thebibliography}

\end{document}